\documentclass[a4paper,10pt]{article}
\usepackage[T1]{fontenc}
\usepackage{lmodern,microtype}
\usepackage[margin=20mm]{geometry}
\usepackage{amsmath,amssymb,amsthm,mathtools,mathrsfs}
\usepackage{array,booktabs,enumitem,needspace,xcolor}
\usepackage[authoryear,square]{natbib}
\setcitestyle{aysep={},citesep={;}}
\usepackage[colorlinks=true,linkcolor=blue!55!black,
  citecolor=blue!55!black,urlcolor=blue!55!black]{hyperref}
\let\cite\citep
\let\citet\citep
\newtheorem{theorem}{Theorem}[section]
\newtheorem{proposition}[theorem]{Proposition}
\newtheorem{lemma}[theorem]{Lemma}
\newtheorem{corollary}[theorem]{Corollary}
\newtheorem{conjecture}[theorem]{Conjecture}
\theoremstyle{definition}
\newtheorem{definition}[theorem]{Definition}
\theoremstyle{remark}
\newtheorem{remark}[theorem]{Remark}
\newcommand{\C}{\mathbb C}
\newcommand{\Q}{\mathbb Q}
\newcommand{\Acon}{A_{\mathrm{con}}}
\newcommand{\Spec}{\operatorname{Spec}}

\newcommand{\End}{\operatorname{End}}

\newcommand{\nrJac}{\operatorname{Jac}}
\newcommand{\nrJdim}{\operatorname{Jdim}}
\newcommand{\nrord}{\operatorname{ord}}
\newcommand{\nrcyc}{\mathrm{cyc}}
\newcommand{\nrGeo}{\operatorname{Geo}}

\numberwithin{equation}{section}
\title{On the Type E case of the Brown--Wemyss realization conjecture}
\author{Shutao Bi\thanks{Email:
  \href{mailto:shutaobi@mail.sdu.edu.cn}{\texttt{shutaobi@mail.sdu.edu.cn}}. Homepage: \url{https://shutaobi.github.io/}.}}
\date{October 2026 (revised November 2026)}
\hypersetup{pdftitle={On the Type E case of the Brown--Wemyss realization conjecture},pdfauthor={Shutao Bi},
  pdfsubject={Brown--Wemyss Numerical Realization, contraction algebras,
  formal potentials and cDV crepant contractions},
  bookmarksnumbered=true,bookmarksopen=false}
\begin{document}
\maketitle
\begin{abstract}
We consider the remaining Type E case of the Brown--Wemyss numerical realization conjecture. The noncommutative splitting lemma and the Type A and D classifications of Brown and Wemyss leave two-generator cyclic potentials with cubic part \(x^3\). For these remaining potentials, we construct geometric realizations under the hypothesis of radical growth dimension at most one. The construction combines marked ADE deformations and cyclic coordinate changes for lengths three, four, and six with a filtered finiteness argument in length five. Under the numerical growth hypothesis, cubic lower bounds exclude the other branches for every choice of formal tail. The resulting contractions have a regular source and a complete local cDV target. Finite-dimensional algebras are realized by threefold flops, and the sole growth-one algebra in this residual case is \(\nrJac(x^3+xy^3)\), for which we give a divisor-to-curve contraction.
\end{abstract}
\clearpage
\begingroup\small\tableofcontents\endgroup
\clearpage
\section{Introduction}

For a flopping curve \(C\) on a threefold, the contraction algebra represents the noncommutative deformation functor of \(\mathcal O_C(-1)\) \cite[Theorem~1.2]{DonovanWemyss}. The realization problem asks which completed associative algebras arise geometrically as contraction algebras of birational crepant contractions. Brown and Wemyss formulated a numerical criterion for realizing completed Jacobi algebras, expressed through the growth of their radical quotients \cite[Conjecture~1.11]{BrownWemyss}. This criterion encompasses both flops, whose contraction algebras are finite-dimensional, and contractions of a divisor to a curve, whose contraction algebras have radical growth dimension one.

The splitting lemma and the Type A and D classifications established by Brown and Wemyss leave two-generator cyclic potentials whose cubic part is \(x^3\) \cite[Corollaries~1.12 and~8.18]{BrownWemyss}. Our contribution concerns this remaining Type E case. The low-degree terms determine the geometric construction. The argument applies to every higher formal tail.

Let \(F_d\) be the free associative algebra on \(x_1,\ldots,x_d\), completed with respect to the ideal generated by these variables. Its augmentation sends every \(x_i\) to zero. We set
\[
 F_d=\C\langle\!\langle x_1,\ldots,x_d\rangle\!\rangle,
 \qquad \mathfrak m_d=(x_1,\ldots,x_d),
 \qquad (F_d)_{\nrcyc}=F_d/\overline{[F_d,F_d]}.
\]
The overline in the cyclic quotient denotes the topological closure of the linear commutator subspace. A cyclic potential is an element of \((F_d)_{\nrcyc}\), represented by a formal series. Its cyclic derivative is defined on a word by
\[
 \partial_{x_i}(x_{j_1}\cdots x_{j_n})
 =\sum_{k:\,j_k=i}x_{j_{k+1}}\cdots x_{j_n}x_{j_1}\cdots x_{j_{k-1}},
\]
and extended linearly and continuously. Thus one deletes each occurrence of \(x_i\) and starts the remaining word immediately after that occurrence. Throughout this paper, every derivative ideal is understood to be a closed two-sided ideal. For a cyclic potential \(\varphi\), we write
\[
 A_\varphi=\nrJac(\varphi)
 =F_d/\overline{(\partial_{x_1}\varphi,\ldots,
                  \partial_{x_d}\varphi)},\qquad
 J_\varphi=\operatorname{rad}(A_\varphi).
\]
For an augmented complete local algebra \(A\) with Jacobson radical \(J\), the radical growth dimension is defined by
\[
 \nrJdim A=\inf\left\{r\geq0:
  \dim_\C(A/J^n)\leq C_r n^r\text{ for every }n\geq1
  \text{ and some }C_r>0\right\}.
\]
When the potential is fixed, we write $A$ for its Jacobi algebra and $J$ for its Jacobson radical. We write \(O_q\) for terms of ordinary word length at least \(q\), and the \(q\)-jet is the truncation through word length \(q\). The order \(\nrord\) of a nonzero series is its least word length, with \(\nrord(0)=\infty\). When positive weights are assigned to the generators, the weight of a word is the sum of its letter weights, and \(O_{\mathrm{wt}}(q)\) denotes terms of weight at least \(q\).

\begin{definition}
\label{def:nr-geometric}
Write \(\nrGeo(A)\) if there is a projective birational crepant
morphism \(\pi:X\to\Spec R\), where \(R\) is a three-dimensional
complete local normal cDV complex algebra, \(X\) is regular,
and the reduced closed fiber is \(\mathbb P^1\), together with a
continuous complete-algebra isomorphism \(A\cong\Acon(\pi)\).
Here cDV means compound Du Val: a general hyperplane section through the closed point of the three-dimensional hypersurface has a Du Val singularity. Regularity is the smooth-source convention in this formal local setting.

We follow Karmazyn's convention for the noncommutative crepant resolution (NCCR) \cite[Theorem~3.3 and Section~3.3]{Karmazyn}. For the basic tilting bundle \(E=\mathcal O_X\oplus N\), let \(e_0\) be the idempotent corresponding to the structure sheaf summand, and set
\[
 \Lambda_\pi=\End_X(E)^{\mathrm{op}},\qquad
 \Acon(\pi)=\Lambda_\pi/\Lambda_\pi e_0\Lambda_\pi.
\]
\end{definition}

\begin{conjecture}[Brown--Wemyss Numerical Realization]
\label{conj:nr-full}
For every \(d\geq1\) and every cyclic formal potential
\(\varphi\in\mathfrak m_d^2\), the condition
\(\nrJdim\nrJac(\varphi)\leq1\) implies
\(\nrGeo(\nrJac(\varphi))\).
\end{conjecture}

A constant term in \(\varphi\) has zero cyclic derivatives. If the linear part of \(\varphi\) is nonzero, one of its cyclic derivatives contains a nonzero constant term and is therefore a unit, which makes the completed Jacobi algebra vanish. We may thus restrict attention to potentials of order at least two.

\Needspace{14\baselineskip}
The Type E arguments below, together with the published splitting and Type A and D realization results, give the following statement.

\begin{theorem}[Geometric Realization]\label{thm:nrf-numerical-realisation}
Let $\varphi$ be a cyclic formal potential in finitely many variables,
with no linear term, and let its nonzero augmented Jacobi algebra satisfy
$\nrJdim\nrJac(\varphi)\le1$. Then $\nrJac(\varphi)$ is isomorphic to
the basic contraction algebra of a projective
birational crepant contraction with regular source, complete local normal
three-dimensional cDV target, and reduced closed fiber $\mathbb P^1$.
The finite-dimensional case is flopping; the growth-one case is
divisor-to-curve. In residual Type E the only growth-one algebra is $\nrJac(x^3+xy^3)$.
\end{theorem}

Here residual Type E refers to the potentials \(x^3+O_4\) left by the published reduction. A marking selects one vertex in an ADE Dynkin diagram; its coefficient in the highest root is the marked length \(\ell\). For the geometric families below, \(\ell\) is also the scheme-theoretic multiplicity of the closed curve fiber, as explained in Proposition~\ref{prop:nr-smooth-arc-geometry}. After the cyclic coordinate changes of Lemma~\ref{lem:nr-xx-elimination}, we may write
\[
 f=x^3+a(xy)^2+bxy^3+cy^4+O_5,
 \qquad \Delta=b^2-8ac.
\]
Proposition~\ref{prop:nr-quartic-orbits} gives five orbits of these four-jets. The two orbits with \(\Delta\ne0\) admit a common length-three construction. The alternating quartic \((xy)^2\) requires a condition on the coefficient of \(y^5\), and the pure quartic \(y^4\) requires further weighted normalization.

For the orbit \(x^3+xy^3\), Theorem~\ref{thm:nr-e6-complete} gives the complete normal form \(x^3+xy^3+p(y)\), with \(p\in y^5\C[[y]]\). The first two rows of the table distinguish whether this series \(p\) vanishes.

For the orbit \(x^3+y^4\), Lemma~\ref{lem:nr-y4-reduction} gives the form
\[
 x^3+y^4+\alpha_{17}xyxy^2
 +\alpha_{20}(xy^2)^2+\alpha_{21}(xy)^3
 +O_{\mathrm{wt}}(24),
 \qquad \operatorname{wt}(x,y)=(4,3),
\]
where the remainder has weight at least \(24\). The subscripts indicate the weights of the three displayed higher terms. The following table gives the case division under \(\nrJdim\nrJac(f)\le1\), together with the locations of the constructions and exclusions. The case division applies to arbitrary formal tails, with the conditions evaluated in the stated normal forms.

\begin{center}
\small
\renewcommand{\arraystretch}{1.16}
\begin{tabular}{@{}>{\raggedright\arraybackslash}p{.26\linewidth}
                    >{\raggedright\arraybackslash}p{.23\linewidth}
                    >{\raggedright\arraybackslash}p{\dimexpr.51\linewidth-4\tabcolsep\relax}@{}}
\toprule
\textbf{Quartic orbit (four-jet)} & \textbf{Further condition} & \textbf{Outcome and proof location}\\
\midrule
\(x^3+xy^3\)
& \(p=0\)
& Marked length \(3\), type \(cE_6\);\newline
divisor-to-curve, \(\nrJdim=1\).\newline
Subsection~\ref{subsec:nr-e6} (Theorem~\ref{thm:nr-e6-complete}).\\
\addlinespace
& \(p\ne0\)
& Marked length \(3\), type \(cE_6\); flop.\newline
Subsection~\ref{subsec:nr-e6} (Theorem~\ref{thm:nr-e6-complete}).\\
\midrule
\(x^3+(xy)^2+y^4\)
& Every higher formal tail.
& Marked length \(3\), type \(cE_6\); flop.\newline
Subsection~\ref{subsec:nr-e6} (Theorem~\ref{thm:nr-e6-complete}).\\
\midrule
\(x^3+(xy)^2\)
& The coefficient of \(y^5\) is zero.
& Excluded: \(\nrJdim\ge3\).\newline
Subsection~\ref{subsec:nrf-five-jet} (Theorem~\ref{thm:nr-weighted-growth}).\\
\addlinespace
& The coefficient of \(y^5\) is nonzero.
& Marked length \(5\), type \(cE_8\); flop.\newline
Subsections~\ref{subsec:nrf-finiteness}--\ref{subsec:nrf-realisation} (Theorem~\ref{nrf:thm:finite} and Proposition~\ref{nrf:prop:geom}).\\
\midrule
\(x^3+y^4\)
& \(\alpha_{17}\ne0\)
& Marked length \(4\), type \(cE_7\); flop.\newline
Subsection~\ref{subsec:nr-transfers} (Theorem~\ref{thm:nr-length-four}).\\
\addlinespace
& \(\alpha_{17}=0\),\newline
\(\alpha_{20}\alpha_{21}\ne0\)
& Marked length \(6\), type \(cE_8\); flop.\newline
Subsection~\ref{subsec:nr-transfers} (Theorem~\ref{thm:nr-length-six}).\\
\addlinespace
& \(\alpha_{17}=0\),\newline
\(\alpha_{20}\alpha_{21}=0\)
& Excluded: \(\nrJdim\ge3\).\newline
Subsections~\ref{sec:gr-quiver-growth}--\ref{sec:gr-sextic} (Theorem~\ref{nrg:thm:boundary}).\\
\midrule
\(x^3\)
& Every higher formal tail.
& Excluded: \(\nrJdim\ge3\).\newline
Subsections~\ref{sec:gr-pbw}--\ref{sec:gr-zero-growth} (Theorem~\ref{nrg:thm:zero}).\\
\bottomrule
\end{tabular}
\end{center}

The divisor-to-curve row realizes \(\nrJac(x^3+xy^3)\), the only growth-one algebra in residual Type E up to complete-algebra isomorphism. All rows marked as flops give finite-dimensional contraction algebras. The generic surface types \(cE_6,cE_7,cE_8\) shown in the table are justified in Section~\ref{sec:gr-assembly}.

The constructions for lengths three, four, and six share a common argument in Subsections~\ref{subsec:nr-geometric-engine}--\ref{subsec:nr-three-directions}. We put the potential into a sum of one-variable series in three linear forms. The Jacobi relations identify the three one-variable derivatives with a common central element, and preparation of the resulting equations gives an arc in the marked ADE deformation space. The derivatives of the constant coefficients establish regularity of the source, while the full endomorphism algebra identifies the contraction quotient.

For length five, Section~\ref{subsec:nrf-rank} first establishes finite-dimensionality. If the Jacobi algebra were infinite-dimensional, its associated graded algebra would have generic rank \(25\) over a central polynomial ring. The determinant of right multiplication on a finite free module would then force this rank to be divisible by \(4\). Finiteness permits polynomial truncation and the application of Wemyss's gluing construction. Section~\ref{sec:gr-assembly} combines these constructions with the growth exclusions to prove the realization theorem. The finite exact calculations used in the proof are described in Section~\ref{sec:gr-certificates}.

\subsection*{Relation to earlier work}
\addcontentsline{toc}{subsection}{Relation to earlier work}

The subsequent classification of contraction algebras by Brown and Wemyss also begins with low-order normal forms \cite[Section~3]{BrownWemyssClass}. Their criterion uses \emph{algebraic length}. They require a central element such that quotienting the Jacobi algebra by the closed ideal it generates gives an algebra isomorphic to the corner \(e\Pi e\), where \(\Pi\) is an ADE preprojective algebra and \(e\) is a primitive idempotent. The marked vertex corresponding to \(e\) determines the algebraic length \cite[Definition~1.1, Corollary~8.14 and Theorem~9.2]{BrownWemyssClass}. The numerical hypothesis concerns only the dimensions of the radical quotients.

For the family \(x^3+xy^3+p(y)\), Pan constructed a smooth flop for every nonzero tail \(p\in y^4\C[[y]]\) and a divisorial contraction at \(p=0\) \cite[Proposition~4.11 and Corollary~4.12]{PanE6}. Brown and Wemyss gave another divisor-to-curve realization of \(\nrJac(x^3+xy^3)\) \cite[Proposition~1.7]{BrownWemyssClass}. The constructions here use the marked deformation spaces of Katz and Morrison \cite[Sections~5--6]{KatzMorrison}, Karmazyn's presentations of the universal flopping algebras \cite[Theorem~1.3 and Examples~4.24--4.27]{Karmazyn}, and Van den Bergh's tilting generators \cite[Propositions~3.2.6--3.2.7]{VanDenBergh}. The quotient of the endomorphism algebra by maps factoring through the structure sheaf summand yields the contraction algebra.

At length five, the geometric input is the fifth family in Wemyss's gluing theorem \cite[Section~1.6 and Theorem~1.7]{WemyssTaming}. Brown and Wemyss identified the deformation algebras of the corresponding curve neighborhoods using free necklace potentials \cite[Corollary~6.4]{BrownWemyssDerived}. The proof of Proposition~\ref{nrf:prop:geom} combines the filtered finiteness argument, polynomial truncation and this gluing construction. The same subsection gives a second route through the later normal form of \cite[Theorem~8.3]{BrownWemyssClass}.

\section{Cyclic normal forms and geometric constructions}
\label{sec:numerical-realisation-routes}

The published reduction of Brown and Wemyss leaves a cyclic cubic \(x^3\) and an arbitrary higher tail. After separating the quartic orbits, we construct the length-three, length-four, and length-six families. The alternating quartic and the growth exclusions are treated in the next two sections.
\subsection{Reduction to two generators}
For a nonzero augmented Jacobi algebra with radical \(J\), the corank is \(\dim_\C J/J^2\). The splitting lemma identifies this number with the number of variables left after the nondegenerate quadratic part has been removed.

\begin{proposition}
\label{prop:nr-full-reduction}
Conjecture~\ref{conj:nr-full} is equivalent to its two-generator
restriction
\[
 \varphi=x^3+O_4,\qquad
 \nrJdim\nrJac(\varphi)\leq1.
\]
Cases of corank at most one and two-generator cases with
\(\dim J^2/J^3=2\) are already geometric in both growth dimensions.
\end{proposition}

\begin{proof}
The noncommutative splitting lemma removes the quadratic variables without altering the isomorphism class of the completed Jacobi algebra \cite[Lemma~4.5 and Appendix~A]{BrownWemyss}. Under the growth hypothesis \(\nrJdim\nrJac(\varphi)\leq1\), the corank bound established by Brown and Wemyss leaves \(\dim J/J^2\) at most two. In corank zero the algebra is \(\C\), and in corank one it is \(\C[[y]]\) or \(\C[[y]]/(y^k)\).

In corank two the cubic part is nonzero. Up to a linear change of variables, this cubic form is \(x^3+y^3\), \(x^2y\), or \(x^3\). The cubic parts \(x^3+y^3\) and \(x^2y\) give \(\dim_\C J^2/J^3=2\), whereas the cubic part \(x^3\) gives \(\dim_\C J^2/J^3=3\). The Type A and D constructions in \cite[Sections~5--6 and Corollary~8.18]{BrownWemyss} realize all the cases covered there in the sense of Definition~\ref{def:nr-geometric}, including their infinite-dimensional forms. Geometric realization is invariant under continuous algebra isomorphism. It therefore remains to treat \(x^3+O_4\). Conversely, the full conjecture applies in particular to these two-variable potentials.
\end{proof}

\subsection{The five quartic orbits}
\label{subsec:nr-quartic}

Two series \(f,g\in F_d\) are \emph{cyclically equivalent} if \(f-g\in\overline{[F_d,F_d]}\), where \([F_d,F_d]\) is the linear span of commutators and the closure is taken in the \(\mathfrak m_d\)-adic topology. In particular, cyclic rotations of a word have the same class; for example, \(xy^2-yxy=[xy,y]\). A necklace is a cyclic equivalence class of words. We write \([h]\) for the cyclic class of an expression \(h\); brackets with two entries, such as \([u,v]\), denote the commutator \(uv-vu\).

Two cyclic potentials represented by \(f,g\in F_d\) are \emph{right-equivalent} if there is a continuous \(\C\)-algebra automorphism \(\Phi:F_d\to F_d\), preserving \(\mathfrak m_d\), such that
\[
 \Phi(f)-g\in\overline{[F_d,F_d]}.
\]
Such an automorphism is an invertible formal change of the noncommuting variables fixing the origin. Thus equality of the transformed potentials is taken in the cyclic quotient \((F_d)_{\nrcyc}\). A change of variables is called unitriangular, or tangent to the identity, when \(\Phi(x_i)-x_i\in\mathfrak m_d^2\) for every generator.

In each cyclic reduction step, a necklace is called good if all of its cyclic shifts avoid the forbidden subwords, including across the cyclic boundary, and bad otherwise.

\begin{lemma}
\label{lem:nr-xx-elimination}
Every \(x^3+O_4\) is right-equivalent, by a continuous
unitriangular automorphism fixing \(y\), to a series whose higher
necklaces contain no adjacent \(xx\), including at the cyclic boundary.
The quartic part of the transformed series is
\begin{equation}
\label{eq:nr-quartic-form}
 a(xy)^2+bxy^3+cy^4.
\end{equation}
\end{lemma}

\begin{proof}
At degree \(n\geq4\), we rotate each necklace containing adjacent letters \(xx\) to the form \(x^2h\), where \(h\) has degree \(n-2\). We write the sum of these terms in degree \(n\) as \([x^2h_n]\). The substitution \(x\mapsto x-h_n/3\), with \(y\) held fixed, changes the cyclic cubic \([x^3]\) by \(-[x^2h_n]\) in this degree. Multiple replacements in the cubic and replacements in terms of degree at least four have degree greater than \(n\). Each elimination step therefore preserves all previous cancellations.

The successive substitutions and their inverses stabilize on every finite word quotient \(F_2/\mathfrak m_2^k\). Their inverse limit yields the required continuous unitriangular automorphism. The three quartic necklaces in \eqref{eq:nr-quartic-form} exhaust those with no adjacent letters \(xx\) in cyclic order. By the cyclic chain rule, each cyclic derivative of the transformed potential lies in the image of the closed derivative ideal. Applying the inverse substitution gives equality of the derivative ideals under transport.
\end{proof}

\begin{proposition}
\label{prop:nr-quartic-orbits}
Put \(\Delta=b^2-8ac\). An invertible linear change of variables preserving the cubic term \(x^3\) has the form
\(x\mapsto\zeta x\), \(y\mapsto tx+vy\), where
\(\zeta^3=1\) and \(v\ne0\). After the change of variables and removal of cyclic \(xx\) terms, the transformed potential has quartic coefficients
\begin{align}
\label{eq:nr-quartic-action}
 a'&=v^2(\zeta^2a+\zeta tb+2t^2c),&
 b'&=v^3(\zeta b+4tc),&c'&=v^4c,\\
 \Delta'&=\zeta^2v^6\Delta.\nonumber
\end{align}
Allowing also a nonzero scalar multiple of the potential, the five
four-jet representatives are
\[
 x^3,\quad x^3+(xy)^2,\quad x^3+xy^3,\quad
 x^3+y^4,\quad x^3+(xy)^2+y^4.
\]
Their strata are respectively the zero quartic,
\(c=b=0,a\ne0\), \(c=0,b\ne0\),
\(c\ne0,\Delta=0\), and \(c\ne0,\Delta\ne0\).
\end{proposition}

\begin{proof}
We expand the transformed quartic and remove all cyclic \(xx\) terms using Lemma~\ref{lem:nr-xx-elimination}. For the contribution of the quartic coefficient \(b\), expand \(\zeta x(tx+vy)^3\) and consider the terms linear in \(t\). Only the middle placement of the new \(x\) factor avoids adjacent letters \(xx\) and contributes to \((xy)^2\), with coefficient \(v^2\zeta tb\). The expansion of \(y^4\), after multiplication by its coefficient, contributes the two words \(xyxy\) and \(yxyx\) to the same cyclic class, giving the total coefficient \(2v^2t^2c\). Together with the direct contributions to the other necklaces, this establishes \eqref{eq:nr-quartic-action} and the stated transformation law \(\Delta'=\zeta^2v^6\Delta\).

The orbit classification follows from these coefficient formulas. For \(c=0,b\ne0\), we choose \(t=-\zeta a/b\) to eliminate \(a'\). For \(c\ne0\), we choose \(t=-\zeta b/(4c)\) to eliminate \(b'\). The remaining alternating coefficient is then proportional to \(a-b^2/(8c)=-\Delta/(8c)\). Rescaling the variables and multiplying the potential by a nonzero scalar yield the five four-jet representatives. Multiplying the potential by a scalar preserves the closed derivative ideal.
\end{proof}

\begin{remark}
\label{rem:nr-quartic-scope}
Each representative specifies only a four-jet. The higher terms in \(\mathfrak m_2^5\) remain arbitrary.
\end{remark}

\subsection{A realization criterion}
\label{subsec:nr-geometric-engine}

Let \(P,Q,R\in u\C[[u]]\), and let \(U,V,W\) be their antiderivatives with zero constant term. The following criterion is applied in each of the length-three, length-four, and length-six constructions.

\begin{theorem}
\label{thm:nr-star-realisation}
Suppose the exact orders of \(P,Q,R\) are one of
\[
 (3,3,2),\qquad(2,4,3),\qquad(2,3,5),
\]
with nonzero leading coefficients.  The Jacobi algebra of
\(U(x)+V(y)+W(-x-y)\) is geometric.  The marked lengths are
respectively three, four and six.  The length-four and length-six contractions are small, meaning that their exceptional loci contain no divisors, and their contraction algebras are finite-dimensional. In length three, the criterion \eqref{eq:nr-e6-divisorial} distinguishes divisor-to-curve contractions from small contractions after the arc coefficients have been constructed.
\end{theorem}

The parameter of the arc will be the common central value of \(P(x),Q(y),R(-x-y)\). We first establish the geometric properties of an arc and the endomorphism comparison. Subsection~\ref{subsec:nr-three-directions} then constructs the arc by one-variable preparation and completes the proof of Theorem~\ref{thm:nr-star-realisation}.

For each length \(\ell\in\{3,4,5,6\}\), let \(H_\ell=\Spec\mathbb H_\ell\) be Karmazyn's marked deformation space, so that \(\mathbb H_\ell\) is its coordinate ring. Write \(\mathcal Y_\ell\to\Spec R_\ell\) for the simultaneous partial resolution over \(H_\ell\), and \(\Lambda_\ell\) for its full endomorphism algebra before quotienting by the structure sheaf projector \cite[Theorem~1.3 and Section~3.3]{Karmazyn}. Let \(e_0\) denote the projector onto the structure sheaf summand in \(\Lambda_\ell\). The universal contraction quotient is
\[
 C_\ell=\Lambda_\ell/\Lambda_\ell e_0\Lambda_\ell.
\]
Take a formal arc \(\Spec\C[[w]]\to H_\ell\) through the marked origin, specified by a ring homomorphism \(\nu:\mathbb H_\ell\to\C[[w]]\). We pull the family back along this morphism, localize and complete its affine target at the central singular point, and base change the projective source to that completed target. Write \(\pi:X\to\Spec R\) for the resulting morphism.

\begin{lemma}
\label{lem:nr-tilting-base-change}
For these marked families, the full endomorphism algebra commutes with
arc base change and target completion.  If the resulting source is
regular and the morphism is crepant with reduced closed fiber
\(\mathbb P^1\), its basic contraction algebra is isomorphic
to the completed \(e_0\)-quotient.
\end{lemma}

\begin{proof}
The comparison requires ordinary endomorphism base change, identification of the basic generator, and passage to the projector quotient.

\noindent\emph{(i) Endomorphisms under base change.}
We write \(\Lambda_\ell=\End(E)^{\mathrm{op}}=\End(E^\vee)\) for Karmazyn's specified tilting bundle \(E=\mathcal O\oplus N\) \cite[Theorem~3.3 and Section~3.3]{Karmazyn}. Both the universal source and the full algebra are flat over \(\mathbb H_\ell\). On a finite affine open cover of the separated source, the \v{C}ech complex for \(\mathcal End(E)\) is bounded with flat terms. Its only nonvanishing cohomology is the flat endomorphism algebra in degree zero. Descending from the last term shows that every boundary module is flat. Tensoring the complex with the arc ring therefore computes the ordinary endomorphisms after pullback and gives zero higher self-Ext. 

Flat base change on the target gives the same comparison after localization and completion. The pulled-back bundle also generates. If a complex on the new source has zero derived Hom from this bundle, adjunction gives zero derived Hom from \(E\) to its pushforward on the original source. Since \(E\) generates there, this pushforward vanishes. Pushforward along the affine base-change map is conservative, so the complex on the new source vanishes as well.

\par\smallskip\noindent\emph{(ii) The generator on the contraction.}
For the pulled-back morphism \(\pi:X\to\Spec R\), assume the regularity and crepancy hypotheses of the lemma. We use \(E\) also for the pulled-back tilting bundle. Rationality of the cDV target gives \(R\pi_*\mathcal O_X=\mathcal O_{\Spec R}\), and the fibers have dimension at most one. The dual bundle \(\mathcal P=E^\vee\) is globally generated, its determinant is relatively ample, and the structure sheaf is a direct summand. We have \(H^1(\mathcal P^\vee)=H^1(E)=0\) because \(E=\mathcal Hom(\mathcal O,E)\) is a summand of \(\mathcal End(E)\). These are the hypotheses of \cite[Propositions~3.2.6--3.2.7]{VanDenBergh}: \(\mathcal P\) is a projective generator of the category of perverse coherent sheaves of perversity minus one for \(\pi\), and \(E\) is a generator of the corresponding category of perversity zero. The endomorphism quotient defining the contraction algebra is taken with the structure sheaf summand distinguished.

It remains to verify the multiplicities of the summands. Modulo \(w\) they reduce to the two distinct indecomposable summands of the basic generator on the completed central partial resolution. Since \(w\) belongs to the Jacobson radical of the full algebra over the complete local target, a nonzero idempotent cannot reduce to zero. A nontrivial direct-sum decomposition upstairs would therefore produce a nontrivial decomposition on the central surface. An isomorphism between the projective summands would likewise descend to that surface. Thus both summands occur with multiplicity one, so \(E\) is basic. Its structure sheaf summand determines the idempotent \(e_0\) used in Definition~\ref{def:nr-geometric}.

\par\smallskip\noindent\emph{(iii) The quotient and its completion.}
For a commutative scalar extension \(H\to S\) and a full \(H\)-algebra \(\Lambda\), the image of the multiplication map \(\Lambda e_0\otimes_H e_0\Lambda\to\Lambda\) is the two-sided ideal \(\Lambda e_0\Lambda\). Forming its cokernel commutes with arbitrary scalar extension, yielding
\[
 (\Lambda/\Lambda e_0\Lambda)\otimes_H S
 \cong(\Lambda\otimes_H S)\big/
        (\Lambda\otimes_H S)e_0(\Lambda\otimes_H S).
\]
For the topology comparison, \(\Lambda\) now denotes the full algebra over the completed target. Let \(\mathfrak n\) be the maximal ideal of that target and \(J\) the Jacobson radical of \(\Lambda\). We have \(\mathfrak n\Lambda\subseteq J\). In the finite-dimensional algebra modulo \(\mathfrak n\), the radical is nilpotent, giving \(J^h\subseteq\mathfrak n\Lambda\) for some \(h\). The two filtrations are therefore cofinal. The ideal \(\Lambda e_0\Lambda\) is a finite module over the completed target ring and is closed for the equivalent \(\mathfrak n\)-adic and \(J\)-adic topologies. Hence the geometric quotient agrees with the augmentation-complete quotient in the loop presentation.
\end{proof}

Deleting the marked vertex from the central ADE diagram leaves chains of type \(A_{m-1}\). These chains describe the residual surface singularities on the central partial resolution. Near each such singularity, the Katz--Morrison projection gives a monic polynomial \(p_m(Z,w)\) for the local equation \(UV=p_m(Z,w)\). The polynomial is centered, meaning that its \(Z^{m-1}\) coefficient is zero, and its lower coefficients vanish at \(w=0\).

\begin{lemma}
\label{lem:nr-smooth-first-jet}
For each residual \(A_{m-1}\) component, write its orthogonal
projection polynomial as
\(p_m(Z,w)=Z^m+\sum_{i=0}^{m-2}a_i(w)Z^i\).
The source \(X\) obtained by pullback and target completion is regular along the central fiber if and only if
\(a_0'(0)\ne0\) at every such component.
\end{lemma}

\begin{proof}
At each residual surface singularity along the central fiber, the orthogonal projection of Katz and Morrison yields the local completed equation \(UV=p_m(Z,w)\) \cite[Section~6, Theorem~3]{KatzMorrison}. The linear part of this defining equation is \(-a_0'(0)w\). When \(a_0'(0)\ne0\), the formal implicit function theorem allows us to eliminate \(w\), leaving the three regular parameters \(U,V,Z\). If \(a_0'(0)=0\), the linear part vanishes identically, and the three-dimensional completed local ring has embedding dimension four, making the point nonregular. 

Away from these residual singularities, the central surface is smooth, and flatness over the regular one-dimensional arc ring ensures regularity at all remaining points. The image of the nonregular locus in the completed local target is closed by projectivity. If nonempty, that image contains the closed point. Regularity along the closed fiber therefore gives regularity on the entire source.
\end{proof}

\begin{proposition}
\label{prop:nr-smooth-arc-geometry}
Assume that, at every residual \(A_{m-1}\) component, the constant coefficient \(a_0(w)\) of the projection polynomial \(p_m(Z,w)\) in Lemma~\ref{lem:nr-smooth-first-jet} satisfies \(a_0'(0)\ne0\).
Then pulling back the marked family along \(\nu\) and completing its target at the central singular point gives a projective birational crepant contraction with regular source, complete local normal cDV hypersurface target and reduced closed fiber \(\mathbb P^1\).  In the primitive length-four, length-five, and length-six families, the contraction is small.  In length three, write the pulled-back quotient parameters as functions \(B_1,B_0,C_1,C_0,D_0,t_K\in\C[[w]]\), so that the centered polynomial relations are
\[
 \beta^3=B_1\beta+B_0,\quad
 \gamma^3=C_1\gamma+C_0,\quad \delta^2=D_0,
 \quad\beta+\gamma+\delta=t_K/3.
\]
With \(B_0',C_0',D_0'\) nonzero at \(w=0\), the contraction \(\pi\) is divisorial precisely when
\begin{equation}
\label{eq:nr-e6-divisorial}
 t_K=0,\qquad B_1=C_1,\qquad B_0=-C_0.
\end{equation}
\end{proposition}

\begin{proof}
Regularity of the source is supplied by Lemma~\ref{lem:nr-smooth-first-jet}. The distinction between the small and divisorial contraction types is governed by whether any positive roots of the ambient ADE root system vanish identically along the formal arc.

\noindent\emph{(i) The target and the closed fiber.}
The target space \(\Spec R\) is flat over \(\C[[w]]\), and its central Cartier section is a Du Val surface singularity. Hence the target is a three-dimensional cDV hypersurface. To verify that \(R\) is an integral domain, divide two nonzero elements by their respective maximal powers of \(w\). Because the central Du Val surface is an integral domain, the reductions of these factored elements have a nonzero product in the central fiber, which guarantees that their product in \(R\) is nonzero. At every smooth point of the central surface, the total space is regular. A codimension-one singular component would have dimension two. Its intersection with the central Cartier section would be positive-dimensional, contradicting the isolated singularity of the Du Val surface. This establishes regularity in codimension one. As a hypersurface, the target automatically satisfies Serre's condition \(S_2\), and is therefore normal.

The simultaneous partial resolution is projective and restricts to an isomorphism over the smooth locus of the surface. Its scheme-theoretic closed fiber is preserved under arc base change and completion of the target. On the central marked surface, the reduction of this fiber is \(\mathbb P^1\), and its scheme-theoretic multiplicity is the coefficient \(\ell\) of the marked simple root in the highest root of the ambient Dynkin diagram \cite[Sections~5--6]{KatzMorrison}. By properness, every connected component of the regular source meets this irreducible closed fiber. Consequently, the source is integral, and the isomorphism over the dense smooth locus establishes birationality of the contraction.

\par\smallskip\noindent\emph{(ii) Roots on the generic surface.}
Over \(K=\C((w))\), we pass to a finite extension on which the deformation parameter lifts to the Cartan space \(\mathfrak h\) of the ADE root system. Each root is a linear form on \(\mathfrak h\), and we evaluate it at this lifted parameter. A root is said to vanish, or to persist along the arc, when this value is zero. The exceptional curves of the generic partial resolution are determined by these vanishing roots. A vanishing root with nonzero coefficient at the marked vertex corresponds to an exceptional curve on the split generic surface. The root calculation below determines when this occurs. In part (iv) we check the closures at the completed target and deduce the contraction type.

For each deleted Dynkin chain of type \(A_{m-1}\), let \(z_1,\ldots,z_m\) be the roots of its centered projection polynomial; thus \(\sum_i z_i=0\). We call this chain an arm of degree \(m\). The root subsystem supported on the deleted diagram is the Levi root system; its integral span is the Levi lattice. The Levi Weyl group is generated by reflections in these roots and permutes the centered root coordinates on each \(A\)-type chain. The first-jet condition \(a_0'(0)\ne0\) ensures that this projection polynomial is Eisenstein in \(w\). Setting \(w=q^m\), substituting \(Z=qz\), and dividing by \(q^m\), the polynomial reduces modulo \(q\) to \(z^m=c\) with \(c\ne0\). By Hensel's lemma, each simple root lifts uniquely to \(\C[[q]]\), the series ring with \(q=w^{1/m}\). The leading terms of the roots \(z_j\) in cyclic order are given by \(c^{1/m}\zeta_m^j w^{1/m}\), where \(\zeta_m\) is a primitive \(m\)-th root of unity. The inertia group of the field extension permutes these \(m\) roots cyclically. In particular, every Levi root is an entry difference \(z_i-z_j\) with \(i\ne j\), which has a nonzero leading term and therefore cannot vanish identically along the arc.

Let \(k>0\) be the coefficient of the marked simple root in an ambient positive root. The evaluation of this root along the arc decomposes as the sum of a base contribution in \(K\) and one contribution \(\sum_j a_jz_j\) from each deleted Dynkin chain. The coefficients \(a_j\) are integers determined up to a common additive integer. If the chain attaches to the marked vertex at node \(h\), these coefficients satisfy the congruence
\[
                 \sum_j a_j\equiv-kh\pmod m.
\]
On this chain, the marked simple root pairs to \(-1\) at the attachment node \(h\) and to zero at all other nodes, so its restriction has weight \(-\omega_h\). On a centered list of roots, the fundamental weight \(\omega_h\) is represented by the sum of \(h\) coordinate entries. Adding a Levi root alters the coefficient vector by a difference of coordinate basis vectors. This yields the congruence, with identical divisibility tests regardless of the orientation of the chain.

Suppose an arm degree \(m\) is prime and coprime to every other arm degree. If an ambient root vanishes identically along the arc, the contribution of this arm must lie in \(K\), which is the intersection of its Puiseux field with the compositum of the other arm fields. Because the root list is centered, the Galois average of this contribution over \(K\) is zero. Hence the contribution vanishes identically. The coefficient of \(w^{1/m}\) in this contribution must therefore satisfy \(\sum_j a_j\zeta_m^j=0\). For prime \(m\), the minimal polynomial of a primitive \(m\)-th root of unity forces all coefficients \(a_j\) to be equal, which implies \(m\mid kh\). For \(m=4\), the vanishing of the leading coefficient forces \(a_1=a_3\) and \(a_2=a_4\), which requires their sum to be even. Applying these divisibility and parity criteria to the marked Dynkin diagrams of \cite[Examples~4.24--4.27]{Karmazyn} yields the conditions in the table below, where \(1\le k\le\ell\).
\[
\begin{array}{ccl}
\toprule
\ell & \text{arm degrees} & \text{necessary conditions for persistence}\\
\midrule
3 &(2,3,3)&2\mid k,\ \text{so }k=2\\
4 &(2,3,4)&3\mid k,\ \text{so }k=3;\ \text{the quartic leading term survives}\\
5 &(4,5)&5\mid k,\ \text{so }k=5;\ \text{the quartic leading term survives}\\
6 &(2,3,5)&2\mid k,\quad3\mid k,\quad5\mid k\\
\bottomrule
\end{array}
\]
The quartic chain in lengths four and five attaches at an end node. An odd coefficient \(k\) yields an odd sum of coefficients, so the leading term of order \(w^{1/4}\) survives. In length four, no other arm produces terms of that fractional order. In length five, orienting the quintic chain with \(h=2\), which is coprime to five, forces its contribution to vanish, leaving the quartic term uncancelled. In length six, no permissible coefficient \(1\le k\le 6\) meets all three divisibility tests simultaneously. This excludes persistent roots in lengths four through six.

\par\smallskip\noindent\emph{(iii) The nine roots in length three.}
In the vertex numbering of \cite[Example~4.24]{Karmazyn}, the marked node \(6\) is adjacent to the nodes \(1,3,5\). To identify the candidate persistent roots, we begin with the root \(2\alpha_6+\alpha_1+\alpha_3+\alpha_5\) belonging to the central \(D_4\) subdiagram. When ambient roots are normalized to have squared norm two, the projections of this root to the two \(A_2\) arms are extremal fundamental weights of squared norm \(2/3\), while its projection to the \(A_1\) arm vanishes.

Every other root with coefficient \(k=2\) differs from this representative by an element of the Levi lattice, so its projection orthogonal to the Levi span is fixed. Within each relevant \(A_2\) weight coset, the three permutations of the fundamental weight are precisely the centered triples of minimal norm. The chosen representative attains this minimum simultaneously in each arm. Any other ambient root of the same total norm must likewise attain the minimum in each arm. Permuting these triples independently yields all nine roots in a single Levi orbit.

Let \(b_i\) and \(c_j\) denote the roots of \(Z^3-B_1Z-B_0\) and \(Z^3-C_1Z-C_0\), respectively. These are the two centered cubics in the statement. In these centered coordinates, the vanishing conditions for the nine roots in the orbit take the form
\[
                    -b_i-c_j=2t_K/3 .
\]
Let \(t_i\) denote the simple-root values in Karmazyn's coordinate system. Substituting \(t_K=t_0+t_1/2\) and \(3t_6=-(t_0+2t_1+t_2+2t_3+t_4+2t_5)\) into the representative yields the expression \(-2t_K/3-(2t_2+t_3)/3-(2t_4+t_5)/3\). Here the last two terms are coordinate entries of the centered cubic lists. Applying Levi permutations generates the remaining eight equations.

Averaging any one of these nine equations over the action of cubic inertia forces \(t_K=0\). Because the cubics are irreducible, the existence of a single opposite pair forces the two root sets to be negatives of each other. Comparing coefficients of the polynomials gives the divisorial conditions \eqref{eq:nr-e6-divisorial}. Conversely, when \eqref{eq:nr-e6-divisorial} holds, we obtain three opposite pairs. The three corresponding roots are mutually orthogonal. For distinct pairs both entries change, and within each fundamental triple the inner product drops by one from the squared norm. The overall pairing between distinct roots drops from two to zero, producing the generic singularity configuration \(A_1^3\).

\par\smallskip\noindent\emph{(iv) The contraction type and crepancy.}
When no ambient root persists, the generic surface is smooth. All exceptional curves of the contraction are then confined to the central fiber, which means that the morphism is small. The target is normal and Gorenstein, and the absence of exceptional prime divisors implies that the contraction is crepant.

In the exceptional length-three case where \eqref{eq:nr-e6-divisorial} holds, the closures of the three generic exceptional curves map onto curve germs in the target after a finite base change \(\Spec\C[[q]]\to\Spec\C[[w]]\) on which the roots split. Here a splitting trait means this spectrum of a discrete valuation ring. To verify that these curve germs survive target completion, consider the relative critical algebra: the coordinate algebra obtained by imposing the ADE surface equation and its partial derivatives in the surface coordinates. Its central fiber is the finite Jacobian algebra of the Du Val singularity. Homogeneous lifts of a basis generate the critical algebra over the polynomial marked base by induction on weight. The unique support point of the central fiber is the vertex singularity. A split critical point therefore has integral coordinates on its trait and specializes to the vertex, ensuring that the morphism extends to the completion.

Each exceptional divisor dominates the parameter curve because the central source surface is integral and maps birationally to its image. On the generic surface, adjunction guarantees that the canonical divisor has zero intersection with each exceptional \((-2)\)-curve. Because the intersection matrix of the configuration is negative definite, all discrepancies vanish over the splitting field and consequently after Galois descent. This establishes crepancy. For the surface-deformation and root descriptions used in this argument, we refer to \cite[Sections~5--6]{KatzMorrison}.
\end{proof}

\begin{remark}
\label{rem:nr-support}
Flatness of the full algebra need not pass to its projector quotient. For example, the \(\C[[w]]\)-free order
\(\left(\begin{smallmatrix}\C[[w]]&\C[[w]]\\
w\C[[w]]&\C[[w]]\end{smallmatrix}\right)\)
has quotient \(\C[[w]]/(w)\) by the two-sided ideal of the diagonal idempotent \(\operatorname{diag}(1,0)\).
Likewise, a nonzero finite local algebra on which a radical central element acts nilpotently cannot be flat over the discrete valuation ring of its parameter.
This is why we verify that endomorphism formation commutes with base change on the full algebra before passing to the contraction quotient.
\end{remark}

\subsection{Constructing an arc in the marked ADE deformation space}
\label{subsec:nr-three-directions}

For \(F=U(x)+V(y)+W(-x-y)\) as in Theorem~\ref{thm:nr-star-realisation}, we seek a formal arc \(\nu:\mathbb H_\ell\to\C[[w]]\) whose completed contraction quotient is isomorphic to \(\nrJac(F)\). Lemma~\ref{lem:nr-common-value} supplies the parameter \(w\) by identifying a common central value of \(P(x),Q(y),R(-x-y)\) in \(\nrJac(F)\). Lemma~\ref{lem:nr-preparation} then converts each equation \(H(u)=w\), for \(H=P,Q,R\), into a monic polynomial equation with no next-to-leading term. The polynomial coefficients and the variable shifts define \(\nu\) in Karmazyn's marked parameter space. The same preparation lemma shows that each constant coefficient has nonzero derivative at \(w=0\), which is the condition in Lemma~\ref{lem:nr-smooth-first-jet} for the pulled-back source to be regular.

\begin{lemma}
\label{lem:nr-common-value}
Let \(P,Q,R\in u\C[[u]]\), and let \(U,V,W\) be their
zero-constant antiderivatives.  The completed Jacobi algebra of
\(F=U(x)+V(y)+W(-x-y)\) is the complete algebra with a central
generator \(w\) and relations
\[
 P(x)=Q(y)=R(z)=w,\qquad x+y+z=0.
\]
\end{lemma}

\begin{proof}
The cyclic derivatives of \(F=U(x)+V(y)+W(-x-y)\) are \(\partial_x F=P(x)-R(-x-y)\) and \(\partial_y F=Q(y)-R(-x-y)\). In the completed Jacobi algebra, setting both cyclic derivatives to zero equates \(P(x)\) and \(Q(y)\) to the common element \(R(-x-y)\). This element commutes with \(x\) because it is equal to the series \(P(x)\), and it commutes with \(y\) because it is equal to the series \(Q(y)\). Because \(x\) and \(y\) generate the algebra topologically, this common value is central. Adjoining the generators \(z=-x-y\) and \(w=P(x)\) yields the displayed presentation. Conversely, eliminating \(z\) and \(w\) via their defining relations defines the inverse homomorphism. Both homomorphisms are continuous with respect to the adic topologies defined by the augmentation ideals.
\end{proof}

\begin{lemma}
\label{lem:nr-preparation}
If \(d\ge2\) and \(H(u)=h_du^d+O(u^{d+1})\), \(h_d\ne0\), then
\(H(u)-w\) is a unit times a unique monic degree-\(d\)
polynomial \(M_H(u,w)\).  Its lower coefficients lie in \(w\C[[w]]\).
After the shift \(v=u+\rho_H(w)\),
where \(\rho_H=m_{d-1}/d\) and \(m_{d-1}\) is the coefficient of \(u^{d-1}\) in \(M_H\), write
\[
 M_H(v-\rho_H,w)=v^d-\sum_{i=0}^{d-2}T_{H,i}(w)v^i.
\]
Then \(T_{H,0}'(0)=1/h_d\).
For \(\widetilde H-H\in u^M\C[[u]]\), \(M\geq d+1\), let \(\widetilde\rho_H\) and \(\widetilde T_{H,i}\) denote the shift and coefficients obtained by applying the same preparation to \(\widetilde H\). Then, for \(0\le i\le d-2\),
\begin{align}
\label{eq:nr-preparation-orders}
 \widetilde\rho_H-\rho_H&\in w^{\lfloor M/d\rfloor}\C[[w]],\\
 \widetilde T_{H,i}-T_{H,i}&\in
 w^{\lceil(M-i)/d\rceil}\C[[w]].\nonumber
\end{align}
\end{lemma}

\begin{proof}
By the formal Weierstrass preparation theorem in \(\C[[w]][[u]]\), the series \(H(u)-w\) factors as an invertible series times a unique monic distinguished polynomial \(M_H(u,w)\in\C[[w]][u]\) of degree \(d\). The roots give both the first jet and the continuity estimates. Set \(w=q^d\). The \(d\) roots of the equation \(H(u)=w\) admit Puiseux expansions \(u_j=\zeta_d^j h_d^{-1/d}q+O(q^2)\) for \(j=0,\ldots,d-1\), where \(\zeta_d\) is a primitive \(d\)-th root of unity. The product \(\prod_{j=0}^{d-1}(u-u_j)\) is invariant under \(q\mapsto\zeta_dq\), so it belongs to \(\C[[w]][u]\). Uniqueness of the distinguished polynomial identifies it with \(M_H(u,w)\). The constant term of \(M_H\) is \(-w/h_d+O(w^2)\). Shifting the variable by \(\rho_H=m_{d-1}/d\) modifies this constant term only by terms in \(O(w^2)\), which yields \(T_{H,0}'(0)=1/h_d\).

To establish the continuity bounds, consider a perturbation of \(H\) by an element of \(u^M\C[[u]]\). We substitute \(u=qv\) and divide the defining equation by \(q^d\). The roots in the rescaled variable \(v\) are simple modulo \(q\), while the perturbation has order at least \(q^{M-d}\). Consequently, each root \(u_j\) changes by an amount of order \(O(q^{M-d+1})\). In the elementary symmetric polynomial giving the coefficient of \(u^i\), the remaining \(d-i-1\) root factors each contribute at least one power of \(q\). The variation in the coefficient of \(u^i\) is therefore divisible by \(q^{M-i}\). Returning to the base parameter \(w=q^d\), this translates into the ceiling bound in \eqref{eq:nr-preparation-orders}. Setting \(i=d-1\) gives the corresponding floor bound for the shift \(\rho_H\). Because \(\rho_H\) and all lower coefficients are divisible by \(w\), the shift preserves these order bounds. A perturbation of the potential of order \(dN+1\) therefore alters the prepared coefficients only within \(w^N\C[[w]]\), for any assignment of complex coefficients.
\end{proof}

\begin{proof}[Proof of Theorem~\ref{thm:nr-star-realisation}]
By Lemma~\ref{lem:nr-common-value}, the completed Jacobi algebra of \(U(x)+V(y)+W(-x-y)\) is presented by the three relations \(P(x)-w=0\), \(Q(y)-w=0\), and \(R(z)-w=0\) along with \(x+y+z=0\), where \(w\) is a central generator. We apply Lemma~\ref{lem:nr-preparation} to each of these three one-variable series and perform the corresponding shifts to eliminate their next-to-leading coefficients. In these centered variables, the linear relation becomes \(\beta+\gamma+\delta=\rho_P+\rho_Q+\rho_R\), so we define the deformation parameter by \(t_K=\ell(\rho_P+\rho_Q+\rho_R)\). The resulting polynomial coefficients match Karmazyn's marked quotient parameters \cite[Section~3.3]{Karmazyn}. Removing the invertible factors preserves the closed derivative ideal, and the inverse coordinate shifts recover the original variables.

Because the leading coefficients of \(P,Q,R\) are nonzero by assumption, Lemma~\ref{lem:nr-preparation} guarantees that \(T_{H,0}'(0)\) is nonzero along each Dynkin arm. The first-jet condition of Lemma~\ref{lem:nr-smooth-first-jet} is therefore satisfied. Proposition~\ref{prop:nr-smooth-arc-geometry} constructs the associated crepant contraction with regular source, and Lemma~\ref{lem:nr-tilting-base-change} identifies its basic contraction algebra with the completed Jacobi algebra. For a small contraction, the contraction algebra is finite over the target ring and supported at the closed point, which proves that it is finite-dimensional over \(\C\). In marked length three, the divisorial alternative is determined by applying criterion \eqref{eq:nr-e6-divisorial} to the prepared coefficients.
\end{proof}

\subsection{A nonzero quartic discriminant}
\label{subsec:nr-e6}

\begin{theorem}
\label{thm:nr-e6-complete}
Let \(f=x^3+O_4\) have reduced quartic
\(a_4(xy)^2+b_4xy^3+c_4y^4\).  If
\(b_4^2-8a_4c_4\ne0\), then \(\nrGeo(\nrJac(f))\).  The construction is a marked length-three
\(E_6\) arc.  If \(c_4\ne0\), it is small. If \(c_4=0\), there exist a continuous \(\C\)-algebra automorphism \(\Phi:F_2\to F_2\) preserving \(\mathfrak m_2\) and a series \(p\in y^5\C[[y]]\) such that
\[
 \Phi(f)-\bigl(x^3+xy^3+p(y)\bigr)\in\overline{[F_2,F_2]}.
\]
The construction is divisorial precisely when \(p=0\).
For every other potential in this stratum, the Jacobi algebra is finite-dimensional.  The sole infinite-dimensional
Jacobi algebra in this stratum, up to complete-algebra isomorphism, is
\(\nrJac(x^3+xy^3)\), of growth dimension one.
\end{theorem}

There are two normal forms, depending on whether \(c_4\) vanishes. The following calculations separate terms that can be removed by coordinate changes from terms that must be supplied by the one-variable summands. For a potential \(g\), the substitution \(x\mapsto x+\varepsilon u\), \(y\mapsto y+\varepsilon v\), with \(\varepsilon^2=0\), changes its cyclic class by \(\varepsilon[u\partial_xg+v\partial_yg]\). We call these first-order variations cyclic tangent terms.

For a graded algebra \(B\), its cyclic quotient, or cocenter, is \(HH_0(B)=B/[B,B]\), where \([B,B]\) is the linear span of commutators. The notation \([b]\) will also denote the image of \(b\in B\) in this quotient. The computations below use polynomial graded algebras before passing to completions.

\begin{lemma}
\label{lem:nr-e6-tangent}
Give \(x,y\) weights \(3,2\), respectively, and set
\(r_x=3x^2+y^3\), \(S=xy^2+yxy+y^2x\).
Every cyclic homogeneous polynomial of weight \(d>9\) has the form
\[
 [u r_x+vS]+\lambda[y^{d/2}],
 \qquad \operatorname{wt}(u)=d-6,\quad
 \operatorname{wt}(v)=d-7,
\]
where the last summand occurs only for even \(d\).
For \(\lambda_0\ne0\), the ordinary graded quotient
\[
 D_{\lambda_0}=\C\langle x,y\rangle/
       (x^2,S+\lambda_0y^3,[x,y^3])
\]
has cyclic quotient spanned by \(1,x,xy\) and pure powers of \(y\).
In degrees at least seven the pure classes vanish unless the degree is
divisible by three.  If \(L_b=y+bx\) and \(t_b=1-\lambda_0b\), then
\begin{equation}
\label{eq:nr-e6-powers}
 [L_b^{3k}]=t_b^k[y^{3k}],\qquad
 [L_b^5]=\frac{1+5t_b^2}{6}[y^5].
\end{equation}
\end{lemma}

\begin{proof}
For the weighted calculation, we identify the cyclic classes remaining modulo the gradient terms \(r_x\) and \(S\). For the ordinary graded calculation, the additional commutator relation \([x,y^3]\) will subsequently be removed by conjugating the generator while fixing the other generator.

\noindent\emph{(i) The quotient by \(r_x,S\).}
In the quotient algebra defined by the relations \(r_x=0\) and \(S=0\), we set \(T=x^2=-y^3/3\). The element \(T\) commutes with both generators \(x\) and \(y\). Over the central polynomial ring \(\C[T]\), we apply the four reduction rules
\[
 xx\mapsto T,\qquad yyy\mapsto-3T,\qquad
 yyx\mapsto-xyy-yxy,\qquad yxyx\mapsto xyxy .
\]
The fourth rule is obtained by comparing \(xS\) with \(Sx\). Each rule decreases the word length or, at equal length, decreases the lexicographic order with \(x<y\). A word is irreducible if it contains none of the left-hand sides. An overlap is a word in which two such left-hand sides intersect; an inclusion occurs when one contains another. These rules have ten proper overlaps and no inclusions. As confirmed by the exact check over \(\mathbb Z[T]\) in Section~\ref{sec:gr-certificates}, the two reductions at each overlap coincide. The irreducible words are
\[
 1,y,y^2,x,xy,xy^2,xyx,xyxy,xyxy^2,yx,yxy,yxy^2.
\]
The overlap checks show that every word reduces to a unique linear combination of these irreducible words. This uniqueness of reduction, called confluence, makes the twelve words a free basis over \(\C[T]\).

We now pass to the cyclic quotient. The relation \([Sy^j]=0\) yields \(3[xy^{j+2}]=0\) for all \(j\ge0\), which implies \(T[x]=T[xy]=0\). The cyclic class of \(xyx\) is \(T[y]\). Similarly, the identities \([Sx]=0\) and \([Sxy]=0\) give \([xyxy]=-2T[y^2]\) and \(2[xyxy^2]=-T[y^3]\). These relations account for all twelve basis words and their central multiples. In weights strictly greater than nine, only pure powers of \(y\) can survive. Rotating each two-sided relation multiple to place the relation factor at the right produces the stated coefficient weights.

\par\smallskip\noindent\emph{(ii) The quotient \(D_{\lambda_0}\).}
In the algebra \(D_{\lambda_0}\), we define \(Z=y^3\) and \(q=xyx\). Because \([x,y^3]=0\), the element \(Z\) is central. Commuting \(x\) with the relation \(S+\lambda_0Z=0\) yields \(x(S+\lambda_0Z)-(S+\lambda_0Z)x=0\). Modulo \(x^2=0\), this simplifies to \(xyxy=yxyx\), which makes \(q\) central with \(q^2=0\). The same twelve words span \(D_{\lambda_0}\) as a module over \(\C[Z]\), acting through the central element \(Z\). Multiplying \(S+\lambda_0Z\) successively on the right by \(1,y,y^2,x,yx\) and taking cyclic equivalence classes yields the following five identities in the cocenter
\begin{align*}
 [xy^2]&=-\lambda_0[Z]/3,&
 Z[x]&=-\lambda_0Z[y]/3,&
 Z[xy]&=-\lambda_0Z[y^2]/3,\\
 q[y]&=(\lambda_0^2/3)Z[y],&
 q[y^2]&=(\lambda_0^2/6)Z[y^2].
\end{align*}
These relations express all twelve cyclic classes in terms of \(1,x,xy\) and pure powers of \(y\). Multiplying the last two relations by \(q\) and using \(q^2=0\) reveals that \(Z^2[y]\) and \(Z^2[y^2]\) vanish. In ordinary degrees seven and higher, the pure powers of \(y\) therefore vanish unless the degree is divisible by three.

To evaluate the powers of the linear form \(L_b=y+bx\), we expand \(L_b^3=t_bZ+b^2q\). Raising this expression to powers and using \([q]=0\) yields the first formula in \eqref{eq:nr-e6-powers}. Multiplying \(L_b^3\) by \(L_b^2\) and substituting the trace identities yields the second formula in degree five. In ordinary degree \(n\geq5\), every cyclic homogeneous polynomial can therefore be written as \([\mu y^n+ux^2+v(S+\lambda_0y^3)+k[x,y^3]]\), where \(u,v,k\) are homogeneous polynomials of degrees \(n-2,n-3,n-4\), respectively.
\end{proof}

\begin{proof}[Proof of Theorem~\ref{thm:nr-e6-complete}]
In this proof, write \(\Delta=b_4^2-8a_4c_4\). The classification divides according to whether the pure quartic coefficient \(c_4\) vanishes. In each case, coordinate transformations absorb all terms lying in the cyclic tangent spaces, while the remaining cyclic classes are matched by one-variable summands in three linear forms, as required by Theorem~\ref{thm:nr-star-realisation}. We then determine which of the resulting algebras are finite-dimensional and identify the unique divisorial member.

\par\smallskip\noindent\textbf{The branch \(c_4=0\).}
When \(c_4=0\), a shear followed by rescaling transforms the reduced quartic into the form \(xy^3\). We then eliminate all cyclic \(xx\) terms by applying Lemma~\ref{lem:nr-xx-elimination}. With respect to the weights \(\operatorname{wt}(x,y)=(3,2)\), all terms beyond the cubic and quartic have weight strictly greater than nine. We apply the first tangent calculation of Lemma~\ref{lem:nr-e6-tangent} in increasing weight. The residual terms are pure powers of \(y\). In weight \(d\), the coordinate adjustments \(x\mapsto x-u\) and \(y\mapsto y-v\) increase the generator weights by \(d-9>0\). Nonlinear contributions produced by multiple replacements have weight at least \(9+2(d-9)>d\), so each cancellation is preserved by subsequent steps. The infinite composition of these substitutions and their inverses converges in the weight-adic topology, bringing the potential to the normal form \(x^3+xy^3+p(y)\).

To construct the geometric realization, set \(\beta=y+x/2\) and \(\gamma=-y+x/2\), and consider the three-direction potential
\[
 G_h=x^3+\beta^4/4-\gamma^4/4+h(\beta),
 \qquad h\in u^5\C[[u]].
\]
Setting \(h=0\) yields the cyclic equivalence class \(x^3+xy^3+x^3y/4\), whose leading error term \(x^3y/4\) has weight eleven. At each even weight \(2m\), Lemma~\ref{lem:nr-e6-tangent} shows that the residual classes are spanned by \([y^m]\). The leading weighted term of \(\beta^m = (y+x/2)^m\) is \(y^m\). Every term containing \(x\) has higher weight. We may therefore choose the coefficient of \(u^m\) in \(h(u)\) to cancel the remaining class \([y^m]\). At odd weights, all cyclic terms lie in the gradient image, so no one-variable correction is required. Induction on weight simultaneously produces the series \(h\) and a convergent automorphism mapping \(G_h\) to the prescribed potential. In terms of the generators \(\beta\) and \(\gamma\), the cubic term is \(x^3 = (\beta+\gamma)^3 = -(-\beta-\gamma)^3\). The three one-variable derivatives therefore have exact orders \((3,3,2)\), so Theorem~\ref{thm:nr-star-realisation} supplies the geometric contraction.

\par\smallskip\noindent\textbf{The branch \(c_4\ne0\).}
When \(c_4\ne0\), we choose a scalar \(\tau\) satisfying the quadratic relation \(a_4+b_4\tau+2c_4\tau^2=0\). The shear \(y\mapsto y+\tau x\) eliminates the alternating quartic term and satisfies \((b_4+4c_4\tau)^2=\Delta\ne0\). Rescaling coordinates and clearing newly created cyclic \(xx\) terms yields \(f_c=x^3+xy^3+cy^4+O_5\) with \(c\ne0\). We then choose a root \(r\) of the quadratic equation
\[
 4cr(1-r)+2r-1=0,
 \quad L_1=y+rx,\quad L_2=-y+(1-r)x,
 \quad a=c(1-r)+1/4,\quad d=cr-1/4.
\]
This quadratic equation excludes the values \(r=0\) and \(r=1\), and determines the coefficients \(a=(1-r)/(4r)\) and \(d=-r/(4(1-r))\). The identities \(a+d=c\), \(4(ar-d(1-r))=1\), and \(ar^2+d(1-r)^2=0\) ensure that the linear combination \(aL_1^4+dL_2^4\) matches the target reduced quartic \(xy^3+cy^4\). All other terms of degree four are removed by a quadratic coordinate adjustment to \(x\).

We then invoke the second calculation of Lemma~\ref{lem:nr-e6-tangent} with \(\lambda_0=4c\) and \(t=1-4cr=r/(1-r)\ne0\). In degrees divisible by three, powers of \(L_1\) supply the residual cyclic class in degree \(3k\) with nonvanishing coefficient \(t^k\). In degree five, \(L_1^5\) and \(L_2^5\) have cyclic classes with coefficients \((1+5t^2)/6\) and \(-(1+5t^{-2})/6\). These coefficients cannot both vanish, because that would require \(t^2=-1/5\) and \(t^2=-5\) simultaneously. Suitable linear combinations in the one-variable series \(h_1(L_1)\) and \(h_2(L_2)\), beginning in degree five, therefore match every surviving cyclic class.

The relation multiples are removed degree by degree. Gradient adjustments eliminate the terms \([ux^2+v(S+4cy^3)]\). To remove the commutator multiple \([k[x,y^3]]\), we apply the exponential of the derivation defined by \(D(x)=[k,x]\) and \(D(y)=0\). This exponential automorphism preserves the cyclic cubic \(x^3\). Here \(n\) is the degree currently being corrected and \(k\) has degree \(n-4\). The exponential is the convergent series \(\exp(D)=\sum_{j\ge0}D^j/j!\): the derivation increases word length by \(n-4>0\). The second variation of the quartic has degree \(2n-4>n\), which enables an inductive cancellation. The coordinate transformations and their inverses converge for arbitrary complex coefficients, producing the three-direction potential
\[
 x^3+aL_1^4+dL_2^4+h_1(L_1)+h_2(L_2).
\]
Because \(L_1+L_2=x\), this potential has derivative orders \((3,3,2)\). Applying Lemma~\ref{lem:nr-preparation} to the derivatives of the two quartic summands gives prepared constant-coefficient derivatives \(1/(4a)\) and \(1/(4d)\). Their sum is \(c/(4ad)\ne0\). Consequently, the length-three divisorial condition \(B_0=-C_0\) fails at first order, which proves that the resulting contraction is small.

\par\smallskip\noindent\textbf{Dimension of the pure-series family.}
Let \(A_p=\nrJac(x^3+xy^3+p(y))\). The element \(T=x^2=-y^3/3\) is central. When \(p=0\), the twelve normal words of Lemma~\ref{lem:nr-e6-tangent} form a free basis for \(A_0\) over \(\C[[T]]\) after completion. The weights of these basis words are bounded and the weight of \(T\) is six, so completing in the weight-adic topology amounts to completing the coefficients over \(\C[T]\). The Jacobson radical modulo \(T\) is nilpotent. Because \(T\in J^2\), the radical filtration and the \(T\)-adic filtration are cofinal with linear bounds on their degrees, which establishes that \(\nrJdim A_0 = 1\).

When \(p\ne0\), the initial relations with respect to the weights \((3,2)\) remain \(3x^2+y^3\) and \(S\). The algebra associated to the filtration by \((3,2)\)-weight is therefore a quotient of the graded algebra defined by \(3x^2+y^3\) and \(S\). Lifting the twelve spanning words and reducing successively in the complete algebra shows that \(A_p\) is generated by at most twelve elements as a module over \(\C[[T]]\).

To prove finiteness, suppose that \(p\ne0\) and the generic algebra \(A_p\otimes_{\C[[T]]}\C((T))\) is nonzero. A simple module over an algebraic closure of \(\C((T))\) has dimension \(r\le3\), because its image in the matrix algebra has dimension \(r^2\le12\). Write \(X,Y\) for the matrices of the generators. Setting \(t^3=T\), the eigenvalues of the diagonalizable matrix \(Y\) are \(\lambda_i=\zeta_i t\), where \(\zeta_i^3=-3\). On the \(\lambda_i\)-eigenspace, the diagonal block of the second cyclic derivative relation reads \(3\lambda_i^2 X_{ii}=-p'(\lambda_i)I\). Hence the trace \(\operatorname{tr}X\) belongs to \(\C((t))\). On the other hand, the relation \(X^2=t^3 I\) implies that \(\operatorname{tr}X\) is an integer multiple of \(t^{3/2}\). Comparing valuations forces \(\operatorname{tr}X=0\), so that the two eigenvalues of \(X\) have equal multiplicity and \(r=2\).

Now expand \(p(y)=h_ny^n+\cdots\) with \(h_n\ne0\) and \(n\geq5\). The two contributions to the trace have leading sum proportional to \(2\zeta_i^{n-3}\) for a repeated eigenvalue of \(Y\), or to \(\zeta_i^{n-3}(1+\omega^{n-3})\) for two distinct eigenvalues, where \(\omega\) is the ratio of the two distinct cube roots and satisfies \(\omega^3=1\). In both cases the leading sum is nonzero, contradicting the condition \(\operatorname{tr}X=0\). The generic algebra must therefore vanish. A sufficiently high power of \(T\) kills all twelve generators, which proves that \(A_p\) is finite-dimensional over \(\C\).

A divisorial contraction in this marked family would have a nonzero contraction algebra at the generic point of its exceptional image curve. On the split generic surface, a generator direct summand pushes forward to a nonfree maximal Cohen--Macaulay module for the contracted \((-2)\)-curve, and its identity endomorphism is nonzero modulo maps factoring through free modules. This quotient is the stable endomorphism algebra. This rules out the divisorial alternative for the finite-dimensional members of this family. The growth-one calculation already excludes a small contraction for \(p=0\).

\par\smallskip\noindent\textbf{The divisorial member.}
We now realize the unperturbed algebra
\[
 \nrJac(x^3+xy^3).
\]
We select the arc in the marked $E_6$ deformation space defined by
\[
 \beta^3=w,\qquad \gamma^3=-w,\qquad \delta^2=-w/3,
 \qquad\beta+\gamma+\delta=0.
\]
The constant coefficients of all three polynomials have nonzero first derivatives at the origin, which guarantees that the source is regular. The two cubics have opposite sets of roots, producing roots that vanish along the arc and have coefficient two at the marked vertex. Proposition~\ref{prop:nr-smooth-arc-geometry} therefore gives a divisor-to-curve contraction. The associated three-direction potential is $-\delta^3+\beta^4/4-\gamma^4/4$. Substituting $x=\beta+\gamma$ and $y=(\beta-\gamma)/2$ yields the cyclic equivalence class
\[
 x^3+xy^3+\tfrac14x^3y.
\]
It remains to verify that eliminating the unwanted term \(x^3y/4\) introduces no pure series in \(y\). Consider the parity involution \(\iota\) defined by \(\iota(x)=-x\) and \(\iota(y)=y\). The leading potential \(x^3+xy^3\) and the discrepancy \(x^3y/4\) are both odd under \(\iota\). At each weight, we apply the odd projector \((1-\iota)/2\) to the tangent decomposition provided by Lemma~\ref{lem:nr-e6-tangent}. Because \(r_x\) is even and \(S\) is odd under \(\iota\), this projector extracts an odd \(x\)-adjustment and an even \(y\)-adjustment, while completely eliminating the pure \(y\)-power summand, which is even. The resulting coordinate substitution commutes with \(\iota\), so that the error term in the next weight remains odd. The weight-adic limit gives the cyclic potential \(x^3+xy^3\) and identifies the completed Jacobi algebra of the arc.

For another divisor-to-curve model of the same algebra, see \cite[Proposition~1.7]{BrownWemyssClass}.
\end{proof}

\subsection{The pure quartic and the two higher-jet conditions}
\label{subsec:nr-transfers}

When the reduced quartic is the pure fourth power \(y^4\), we assign weights \(\operatorname{wt}(x,y)=(4,3)\) so that both leading terms \(x^3\) and \(y^4\) have weighted degree twelve. Among higher terms of weight below twenty-four, the surviving necklaces have weights \(17,20,21\). Their coefficients govern the realization for marked lengths four and six.

\begin{lemma}
\label{lem:nr-y4-reduction}
Every formal potential \(x^3+y^4+O_5\) is right-equivalent to
\[
 x^3+y^4+\alpha_{17}xyxy^2
       +\alpha_{20}(xy^2)^2+\alpha_{21}(xy)^3
       +O_{\mathrm{wt}}(24),
 \qquad\operatorname{wt}(x,y)=(4,3),
\]
where every remaining necklace avoids the cyclic subwords \(xx\) and \(yyy\).
The corresponding ordinary remainder has order at least seven.
\end{lemma}
\begin{proof}
In weighted degree \(n>12\), rotate each bad necklace to place a forbidden subword at the beginning. The sum of all such bad terms can then be expressed cyclically as \(x^2u_n\) plus \(y^3v_n\). The elementary substitutions \(x\mapsto x-u_n/3\) and \(y\mapsto y-v_n/4\) raise the weights of the generators by \(n-12>0\). They cancel these terms to first order and perturb the potential nonlinearly only in weights strictly greater than \(n\). The substitutions and their inverses therefore converge in the weight-adic topology.

Because the formal tail begins in ordinary degree five, the steps at weights thirteen and fourteen are trivial. In particular, they introduce no \((xy)^2\) term. Every surviving higher necklace must avoid adjacent letters \(xx\) and subwords \(yyy\), so each occurrence of \(x\) must be followed by either one or two letters \(y\). Surviving necklaces are therefore concatenations of the two basic blocks \(xy\) of weight seven and \(xy^2\) of weight ten, containing at least two occurrences of \(x\). Listing all block words of weighted degree below twenty-four yields the three displayed terms of weights seventeen, twenty, and twenty-one. Every remaining higher necklace has weighted degree at least twenty-four and ordinary word length at least seven.
\end{proof}

\begin{theorem}
\label{thm:nr-length-four}
For every \(a\ne0\) and every cyclic formal tail \(h\in O_6\),
\[
 f=x^3+y^4+a\,xyxy^2+h
\]
is right-equivalent to \(x^3+V(y)+W(y-x)\), where
\(V(u)=u^4+O(u^5)\), \(W(u)=(a/5)u^5+O(u^6)\).
Its Jacobi algebra is finite-dimensional and geometric through a smooth
marked length-four small contraction.
\end{theorem}

\begin{theorem}
\label{thm:nr-length-six}
For \(b,c\ne0\), every
\[
 f=x^3+y^4+b(xy^2)^2+c(xy)^3+O_7
\]
is right-equivalent to \(x^3+V(y)+W(y-\lambda x)\), where
\[
 \lambda=-3c/(2b),\quad d_0=b/(3\lambda^2),\quad
 V(u)=u^4+O(u^5),\quad W(u)=d_0u^6+O(u^7).
\]
It is realized by a smooth marked length-six small contraction and has
finite-dimensional Jacobi algebra.
\end{theorem}

The proofs of Theorem~\ref{thm:nr-length-four} and Theorem~\ref{thm:nr-length-six} vary the one-variable series together with a coordinate change. The following lemma gives the order estimates used in these degree-by-degree constructions.

\Needspace{17\baselineskip}
\begin{lemma}
\label{lem:nr-flow-transfer}
Let \(G=U(x)+V(y)+W(L)\), \(L=y-\lambda x\), \(\lambda\ne0\).
Suppose \(U'\) has exact order \(d\geq2\) and
\(W'(u)=bu^e+O(u^{e+1})\), \(e>d\), \(b\ne0\).
For \(k\geq1,s\geq0\) with \(r=d(k-1)+s\geq2\) and for every
\(\rho\in\C\), there exist a modified \(U_\rho\) and an
automorphism \(\Psi_\rho\) tangent to the identity and fixing \(y\), such that
\begin{align*}
 \nrord(U_\rho-U)&\geq dk+s,\\
 \nrord(\Psi_\rho(x)-x),\ \nrord(\Psi_\rho^{-1}(x)-x)&\geq r,\\
 [\Psi_\rho(U_\rho(x)+V(y)+W(L))]
 &=[G]+\rho b^k[L^{ek}x^s]+O_{ek+s+1}.
\end{align*}
The analogous construction in the variable \(y\) assumes \(V'\) has order \(d\), fixes \(x\), and
adds \(\rho b^k[L^{ek}y^s]\) in degree \(ek+s\).
\end{lemma}
\begin{proof}
We vary the one-variable series \(U\) along a one-parameter flow and compensate with a coordinate substitution in \(x\) while holding \(y\) fixed. The noncommutative telescoping identity below ensures that the net variation of the potential becomes a multiple of a power of \(W'(L)\).

Introduce a formal time parameter \(\tau\), set \(U_0=U\), and define the family of potentials \(G_\tau=U_\tau(x)+V(y)+W(L)\). We set \(A_\tau=U_\tau'(x)/\lambda\) and \(B=W'(L)\), and let \(D_\tau\) be the derivation of the completed free algebra determined by its values on \(x,y\) below. We solve the differential system
\begin{align*}
 \partial_\tau U_\tau(u)&=\rho u^s(U_\tau'(u)/\lambda)^k,\\
 C_\tau&=\sum_{i=0}^{k-1}B^{k-1-i}x^sA_\tau^i,\qquad
 D_\tau(x)=-\rho C_\tau/\lambda,\quad D_\tau(y)=0.
\end{align*}

For a solution of this system, let \(\Psi_\tau\) be the coordinate flow satisfying \(\partial_\tau\Psi_\tau=\Psi_\tau D_\tau\) with initial condition \(\Psi_0=\mathrm{id}\). Noncommutative telescoping and cyclic rotation of words yield
\[
 [x^sA_\tau^k]-[B^kx^s]=[C_\tau(A_\tau-B)],
 \quad \partial_\tau[\Psi_\tau(G_\tau)]
           =\rho[\Psi_\tau(B^kx^s)].
\]

The order estimates also give existence of this flow. Because \(r=d(k-1)+s\geq2\), each step in the coefficient recursion for \(\partial_\tau U_\tau\) increases the word order by at least \(r-1>0\). The coefficients of \(U_\tau\) are therefore polynomials in the parameter \(\tau\). The element \(C_\tau\) has order at least \(r\), so the flow equation \(\partial_\tau\Psi_\tau=\Psi_\tau D_\tau\) with \(\Psi_0=\mathrm{id}\) can be solved order by order, yielding coefficients for \(\Psi_\tau\) that are polynomial in \(\tau\). On each finite quotient \(F_2/\mathfrak m_2^N\), the substitution acts unipotently. Its inverse is likewise polynomial in \(\tau\) and satisfies the same order bound. We may therefore evaluate the flow at \(\tau=1\) compatibly on all finite quotients.

The degree-\(ek+s\) term on the right-hand side of the flow equation equals \(\rho b^k[L^{ek}x^s]\) identically throughout the flow. Integrating these polynomial coefficients from \(\tau=0\) to \(\tau=1\) yields the asserted leading term. For the \(y\)-branch, we set \(A_\tau=-V_\tau'(y)\), \(\partial_\tau V_\tau=\rho u^s(-V_\tau')^k\), \(C_\tau=\sum_{i=0}^{k-1} B^{k-1-i}y^sA_\tau^i\), and \(D_\tau(y)=\rho C_\tau\), while leaving \(x\) fixed. The exact flow identity then gives the corresponding leading variation and order estimates.
\end{proof}

\begin{proof}[Proof of Theorem~\ref{thm:nr-length-four}]
We choose the higher-order coefficients of the one-variable series \(V\) and \(W\) successively to match the potential, and then absorb the relation multiples through coordinate changes. The auxiliary quotient below identifies the cyclic classes that must be supplied by \(V(y)\) and \(W(y-x)\). Throughout this proof, principal gradient corrections mean coordinate changes using the derivatives \(3x^2\) and \(4y^3\) of \(x^3+y^4\).

\noindent\textbf{Residual cyclic classes.} Write \(u=xy\), \(v=xy^2\), \(g=xyxy^2\), and \(q_{17}=[x,\partial_xg]\). To determine the degree-\(m\) terms that survive the principal gradient corrections and coordinate adjustments, consider the graded algebra
\[
 D_{17}=\C\langle x,y\rangle/(x^2,y^3,q_{17}).
\]
Every non-forbidden mixed necklace is a cyclic word in the basic blocks \(u\) and \(v\). For \(m\ge6\), we show that the residual cyclic classes are spanned by \([u^{m/2}]\) when \(m\) is even, together with \([u^{2k-2}v^2]\) when \(m=4k+2\). All odd degrees in this range vanish in the cyclic quotient. Modulo the quadratic and cubic monomial relations \(x^2=0\) and \(y^3=0\), direct cyclic differentiation of \(g\) yields
\[
 q_{17}=uv+vu-yv x-y^2u x .
\]
Every nonempty word \(H\) in the blocks \(u\) and \(v\) begins with the letter \(x\). Consequently, the cyclic classes \([Hyvx]\) and \([Hy^2ux]\) vanish because adjacent letters \(xx\) appear across the cyclic boundary. This leaves the cyclic relation \([Huv]=-[Hvu]\), which demonstrates that in any cyclic word containing at least three blocks, exchanging two adjacent unlike blocks changes the sign of the cyclic class. Modulo the monomial relations, we also have \(q_{17}y^2=-yv^2-y^2uv\). Multiplying on the left first by \(y\) and then by \(x\) shows that \(v^3=0\) in \(D_{17}\). Sorting the blocks into the standard form \(u^i v^j\) therefore kills any cyclic class containing three or more occurrences of the block \(v\).

For a sorted mixed word \(u^iv^j\) with at least three blocks, cyclically rotating one block \(u\) from the beginning to the end and commuting it back past the \(j\) factors of \(v\) introduces the sign \((-1)^j\). Moving a factor of \(v\) in the same cyclic manner introduces \((-1)^i\). Hence both exponents \(i\) and \(j\) must be even for the cyclic class to be nonzero. The only mixed two-block word is \(uv\), which has degree five. This establishes the asserted list of residual classes in every degree \(m\ge6\). Rotating the two-sided relation multiples expresses the remaining discrepancy in the form \([Hq_{17}+x^2u_m+y^3v_m]\), where the coefficient polynomials have degrees \(m-5,m-2,m-3\), respectively.

\par\smallskip\noindent\textbf{Matching degree six.}
In the initial potential \(x^3+y^4+(a/5)(y-x)^5\), the unique non-forbidden fifth necklace has coefficient \(a\). We remove all forbidden fifth necklaces using principal corrections in \(x\) and \(y\). Set \(L=y-x\). To adjust the sixth jet, we apply the \(y\)-flow of Lemma~\ref{lem:nr-flow-transfer} prior to those normalizing substitutions, using parameters \((d,e,k,s)=(3,4,1,2)\) and amplitude \(-A/4\). The leading variations are \(Ay^5\) in the series \(V\) and \(-Ay^2/4\) in the coordinate \(y\). By the exact flow identity, their net variation in degree six is \(-aA[L^4y^2]/4\). Because the normalizing substitutions are tangent to the identity, they preserve this leading term.

Expanding \(L^4y^2\) in terms of blocks, the only necklace that avoids forbidden subwords cyclically is \(xy^2xy^2 = v^2\). For the second independent direction, we add \(B\xi^6\) to the one-variable series \(W(\xi)\). This adds \(B[L^6]\) to the potential, which expands to three cyclic rotations of \(v^2\) and two cyclic rotations of \(u^3\) with opposite signs. Combining these two contributions yields
\begin{equation}
\label{eq:nr-length-four-six}
 (-aA/4+3B)v^2-2Bu^3
\end{equation}
in the cyclic quotient. The linear transformation sending \((A,B)\) to the coefficients of \(v^2\) and \(u^3\) has determinant \(a/2\). Since \(a\ne0\), we can match both residual coefficients in degree six.

\par\smallskip\noindent\textbf{Higher degrees and convergence.}
Assume that the current family agrees with the prescribed potential below degree \(m\ge7\). At degree \(m=4k\), adding \(\mu\xi^m\) to \(W(\xi)\) contributes \(2\mu[u^{2k}]\). At degree \(m=4k+2\), the same variation produces the coefficient \(-2\mu\) on \([u^{2k+1}]\). The second independent direction is provided by the \(y\)-flow with exponent \(k\) and shift \(s=2\), whose leading contribution is \(\rho a^k[L^{4k}y^2]\).

The appended factor \(y^2\) forces the final block of the word to be \(v\) after cyclic rotation. For a necklace containing exactly two factors of \(v\), the other pair of consecutive \(y\)'s may appear in any of the first \(2k-1\) positions. Commuting the two factors of \(v\) together produces alternating signs that sum to one. Necklaces containing three or more factors of \(v\) vanish after sorting, by \(v^3=0\), while the presence of the appended pair excludes the pure \(u\)-class. The flow therefore contributes \(\rho a^k[u^{2k-2}v^2]\) to the cyclic quotient. In odd degrees, there are no surviving cyclic classes to match. At each step, we apply the transfer flow to the current family of potentials. The family update has order at least \(3k+2\), while the coordinate map and its inverse differ from the identity by terms of order at least \(3k-1\).

To eliminate the relation multiple \([Hq_{17}]\), we apply the continuous automorphism that conjugates the generator \(x\) by \(\exp(H/a)\) while holding \(y\) fixed. This coordinate substitution preserves the cyclic cubic \(x^3\) and leaves the pure quartic \(y^4\) invariant. The first variation of \(ag\) under this substitution is \([Hq_{17}]\), while the second variation has degree \(2m-5>m\). Principal gradient corrections in \(x\) and \(y\) remove the multiples of \(x^2\) and \(y^3\), with signs chosen to cancel the remaining error terms. The updates to the one-variable series, the coordinate transformations, and their inverses all stabilize on every finite word quotient \(F_2/\mathfrak m_2^N\). Passing to the inverse limit and setting \(X=x\) and \(Y=y-x\) yields a three-direction potential with derivative orders \((2,4,3)\). Theorem~\ref{thm:nr-star-realisation} therefore establishes the claimed smooth marked length-four small contraction.
\end{proof}

\begin{proof}[Proof of Theorem~\ref{thm:nr-length-six}]
We use the same division into residual cyclic classes and relation multiples as in the length-four proof. In this setting, the existence of a central element of degree five reduces the cyclic calculation to a finite-dimensional auxiliary algebra.

\noindent\emph{(i) A finite quotient.}
Set \(L=y-\lambda x\), \(g_6=b(xy^2)^2+c(xy)^3\), \(q_6=[x,\partial_xg_6]\), and \(D_6=\C\langle x,y\rangle/(x^2,y^3,q_6)\). To determine the cyclic classes that must be matched in each degree, we first identify a central element of degree five. In the quotient algebra by the ideal \((x^2,y^3)\), the linear operation \(h\mapsto[x,\partial_xh]\) kills every bad necklace containing forbidden subwords. For a necklace \(x^2a\), cyclic differentiation gives \([x,xa+ax]=[x^2,a]\), while all other terms contain a factor of \(x^2\). For a representative \(y^3a\), every term of its cyclic \(x\)-derivative contains a factor of \(y^3\). The non-forbidden coefficients of \(d_0L^6\) are \(3d_0\lambda^2=b\) and \(-2d_0\lambda^3=c\). Consequently, we have the congruence
\[
 q_6=-6\lambda d_0[x,L^5]\pmod{(x^2,y^3)}.
\]
Because \([y,L^5]=\lambda[x,L^5]\), the element \(z=L^5\) commutes with both generators and is therefore central in \(D_6\).

After rescaling the coordinate \(x\) by \(\lambda\), the quotient algebra by the ideal generated by \(z\) becomes
\[
 B_6=\C\langle x,y\rangle/(x^2,y^3,(y-x)^5).
\]
In this finite calculation, \(x\) denotes the rescaled variable and \(L=y-x\). The resulting spanning classes will be transported back to \(D_6\) by the inverse rescaling. We index the columns in degree \(n\) by words avoiding adjacent letters \(xx\) and subwords \(yyy\). For each pair of non-forbidden words \(a,b\) with \(|a|+|b|=n-5\), we expand the product \(a(y-x)^5b\) and discard all words containing forbidden subwords. The resulting rows span the relations in degree \(n\). For the cyclic cocenter \(HH_0(B_6)\), we identify cyclic permutations of words and set any class with a forbidden cyclic rotation to zero. The rows \([(y-x)^5b]\) then suffice to span the relations, because a left multiplier rotates cyclically to the right. Gaussian elimination over the field of rational numbers gives the vector space dimensions in degrees zero through fifteen
\[
\begin{aligned}
 \dim(B_6)_n &: 1,2,3,4,5,6,6,6,6,6,5,4,3,2,1,0,\\
 \dim HH_0(B_6)_n &: 1,2,2,1,1,0,0,0,0,0,0,0,0,0,0,0 .
\end{aligned}
\]
Adding the rows corresponding to \(1,y,y^2,L,L^2,L^3,L^4\) to the cyclic relation matrices in their respective degrees yields full column rank. Thus the images of these seven elements span \(HH_0(B_6)\). The computer algebra verification in Section~\ref{sec:gr-certificates} confirms these matrix ranks. Because the degree-fifteen component of the algebra \(B_6\) vanishes identically, every word of degree fifteen and higher vanishes as well.

\par\smallskip\noindent\emph{(ii) Lifting the cyclic generators.}
Because \(z\) is central in \(D_6\), we have the isomorphism \(HH_0(D_6)/zHH_0(D_6)=HH_0(B_6)\). For any homogeneous cyclic class in \(D_6\), we subtract a linear combination of the seven lifts of the spanning classes of \(HH_0(B_6)\). The resulting remainder belongs to \(zHH_0(D_6)\), and can therefore be written as \(z\) times a cyclic class of degree lower by five. By induction on degree, these seven elements generate the cyclic cocenter \(HH_0(D_6)\) as a module over the polynomial ring \(\C[z]\).

Consequently, in every degree \(m\ge7\), the residual cyclic classes are spanned by \([L^m]\) and, when \(m=5k+j\) with \(j\in\{1,2\}\), by \([L^{5k}y^j]\). After these surviving classes have been matched, the remaining discrepancy in degree \(m\) takes the form \([x^2u_m+y^3v_m+H_mq_6]\), where the coefficient polynomials have degrees \(m-2,m-3,m-6\), respectively.

\par\smallskip\noindent\emph{(iii) Matching the formal tail.}
The initial potential \(x^3+y^4+d_0L^6\) has the required sextic coefficients \(3d_0\lambda^2\) and \(-2d_0\lambda^3\). We remove its forbidden sixth-degree necklaces using principal corrections in \(x\) and \(y\). Assume by induction that the potential matches the prescribed series through degree \(m-1\). Adjusting the degree-\(m\) coefficient in the one-variable series \(W\) supplies the required multiple of \([L^m]\). To supply the class \([L^{5k}y^j]\), we apply the second transfer flow of Lemma~\ref{lem:nr-flow-transfer} with parameters \(d=3,e=5,s=j\), setting the flow amplitude \(\rho\) equal to the target coefficient divided by \((6d_0)^k\).

The order condition \(r=3(k-1)+j\) is at least two in all relevant degrees. The minimal value of \(k\) is two when \(j=1\) (degree eleven) and one when \(j=2\) (degree seven). The exact flow identity of Lemma~\ref{lem:nr-flow-transfer} therefore applies from the very first step, modifying the series \(V\) only in orders at least \(3k+j\).

Finally, we remove the relation multiple \([H_mq_6]\) by applying the continuous automorphism that conjugates the generator \(x\) by \(\exp(H_m)\) while holding \(y\) fixed. The resulting first variation provides \([H_mq_6]\), while the second variation on \(g_6\) has degree \(2m-6>m\). Principal gradient corrections in \(x\) and \(y\) remove the other relation multiples. Each previously matched degree is preserved throughout these adjustments. As \(m\) tends to infinity, the orders of the family updates and of the differences of the coordinate maps and inverse maps from the identity tend to infinity. Their limits give a continuous right-equivalence in the word-adic topology. Setting \(X=\lambda x\) and \(Y=-y\) yields derivative orders \((2,3,5)\). Theorem~\ref{thm:nr-star-realisation} applies to these orders and constructs the claimed smooth marked length-six small contraction.
\end{proof}

\section{Realization of the length-five stratum}
\label{subsec:nrf-rank}

For the alternating quartic, the growth bound \(\nrJdim\leq1\) forces the pure fifth coefficient to be nonzero. We first normalize this coefficient, and then prove that every Jacobi algebra with the resulting five-jet is finite-dimensional. The geometric realization follows by polynomial truncation and gluing.

\subsection{The alternating quartic and its fifth jet}
\label{subsec:nrf-five-jet}
\begin{theorem}
\label{thm:nr-weighted-growth}
If every word of \(f\) has weight at least six for
\(\operatorname{wt}(x,y)=(2,1)\), then for every \(N\geq1\)
\begin{equation}
\label{eq:nr-cubic-growth}
 \dim_\C\bigl(\nrJac(f)/J^N\bigr)
 \geq\#\{(a,b,c)\in\mathbb Z_{\geq0}^3:2a+b+3c<N\}
 \geq\bigl(\lfloor(N-1)/6\rfloor+1\bigr)^3.
\end{equation}
In particular \(\nrJdim\nrJac(f)\geq3\).
\end{theorem}

\begin{proof}
We first count the normal words in a formal deformation of a fixed homogeneous potential. Specialization of the finite relation matrix then establishes the dimension lower bound at every point in the coefficient space.

Let \(K_N\) denote the closed subspace spanned by words of weight at least \(N\). Because cyclic differentiation reduces \((2,1)\)-weight by at most two, only coefficients of the potential in weights six through \(N+1\) can affect the truncated quotient. We vary these coefficients independently about the homogeneous potential \(f_0=(xy-yx)^2/2\), introducing a formal perturbation \(g=\sum e_i h_i\) over \(T=\C[[e_1,\ldots,e_r]]\). Here the \(h_i\) form a basis of the cyclic words in weights six through \(N+1\), and the \(e_i\) are independent central parameters. We freely adjoin a generator \(\zeta\) of weight three. The cyclic derivatives of the perturbed potential \(\zeta(xy-yx)-\zeta^2/2+g\) produce the reduction rules
\[
 yx\longmapsto xy-\zeta,\quad
 \zeta y\longmapsto y\zeta+\partial_xg,\quad
 \zeta x\longmapsto x\zeta-\partial_yg.
\]
For the fixed weight cutoff \(N\), we also work modulo a power of the parameter ideal \((e_1,\ldots,e_r)\). We will let this parameter cutoff tend to infinity while keeping \(N\) fixed. Each cyclic derivative correction increases the parameter degree. At any fixed parameter degree, the remaining terms strictly decrease the word length or, at equal length, reduce the number of inversions with respect to the ordering \(x<y<\zeta\). This guarantees termination, and none of the reduction rules decreases the weighted degree. The unique overlap among the leading words is \(\zeta yx\). The two distinct reduction paths for this overlap differ by \(-[y,\partial_yg]-[x,\partial_xg]\), which vanishes identically by the cyclic commutator identity \([x,\partial_x g] + [y,\partial_y g] = 0\). Disjoint reductions commute.

The terminating confluent reductions show that the ordered monomials \(x^ay^b\zeta^c\) satisfying \(2a+b+3c<N\) yield a unique normal form at every parameter cutoff. Passing to the inverse limit in parameter degree shows that the truncated algebra is a finite free \(T\)-module with these ordered monomials as a basis. Eliminating the auxiliary generator \(\zeta=xy-yx\) identifies this module with the weighted quotient of the deformed completed Jacobi algebra.

We now extend this dimension bound from the formal neighborhood of \(f_0\) to arbitrary complex coefficients. Consider the finite relation matrix defining the truncated word quotient. Its entries are polynomials in the coefficients \(e_i\). The generic rank of this matrix over the quotient field of \(T\) equals the total dimension of the word space minus the number of ordered monomials. Every minor of size strictly greater than this generic rank vanishes identically as a formal power series at the origin, and because its entries are polynomials in the coefficients, it must vanish identically as a polynomial everywhere. Specialization to any complex coefficient values can only decrease the rank of the relation matrix, which can only increase the dimension of the resulting quotient. Translating the coefficient coordinates by the coefficients of \(f_0\) therefore preserves the same dimension lower bound for every potential of weight at least six. At any fixed weight cutoff, the closed derivative ideal and the ordinary derivative ideal have identical images on the finite-dimensional word space, so the higher-order formal tail introduces no additional relations at that cutoff.

Because the powers of the Jacobson radical satisfy \(J^N\subseteq K_N\nrJac(f)\), the radical quotient surjects canonically onto the weighted quotient. This establishes the first inequality in \eqref{eq:nr-cubic-growth}. For the second inequality, restricting each of the three nonnegative integer exponents \(a,b,c\) to the range \(0\le a,b,c \le \lfloor(N-1)/6\rfloor\) yields the weight bound \(2a+b+3c\le 6\lfloor(N-1)/6\rfloor < N\). The number of such triples is at least \((\lfloor(N-1)/6\rfloor+1)^3\), which completes the proof.
\end{proof}

\begin{corollary}
\label{cor:nr-five-jet}
If \(f=x^3/3+(xy)^2/2+O_5\) and \(\nrJdim\nrJac(f)\leq1\),
its pure \(y^5\) coefficient is nonzero.  A formal coordinate change
and a nonzero scalar multiple put its five-jet into
\[
 x^3/3+(xy)^2/2-xy^4-y^5/5.
\]
\end{corollary}

\begin{proof}
The only cyclic necklaces of degree five avoiding adjacent letters \(xx\) are \(xyxy^2\), \(xy^4\), and \(y^5\). If the coefficient of \(y^5\) were zero, every surviving word in the potential would have \((2,1)\)-weight at least six, which by Theorem~\ref{thm:nr-weighted-growth} would force \(\nrJdim\nrJac(f)\geq3\), contradicting the hypothesis \(\nrJdim\leq1\).

Let \(A,B,C\) denote the coefficients of \(xyxy^2\), \(xy^4\), and \(y^5\), respectively, with \(C\ne0\). Rescaling coordinates by \(x\mapsto s^2x\) and \(y\mapsto sy\), and multiplying the overall potential by \(s^{-6}\) with \(s=-5C\), normalizes the leading pure power coefficient to \(C=-1/5\). We then apply the shear \(y\mapsto y+dx\) and restore the standard four-jet via the quadratic substitution \(x\mapsto x-dxy-d^2x^2/2\). This combined shear and restoration changes the coefficient \(B\) of the necklace \(xy^4\) to \(B-d\), while the restoring substitution affects only necklaces containing at least two occurrences of \(x\). Choosing \(d=B+1\) normalizes this coefficient to \(-1\). We then eliminate the bad fifth necklaces using a cubic coordinate correction in \(x\). Let \(A'\) denote the coefficient of \(xyxy^2\) after these steps. The quadratic transformation \(y\mapsto y-A'y^2\) cancels this remaining term. All remaining terms have ordinary word length at least six.
\end{proof}

\subsection{The length-five finiteness theorem}
\label{subsec:nrf-finiteness}
For the finiteness calculation, we assign weights \(1,2\) to \(x,y\), respectively. We write \(F^nA\) for the image of series of weight at least \(n\), and \(\operatorname{gr}A=\bigoplus_{n\ge0}F^nA/F^{n+1}A\). The initial form of an element is its first nonzero class in this graded algebra. Proposition~\ref{nrf:prop:geom} relates these generators to the preceding five-jet through an explicit change of coordinates. Let
\[
p(x)=a_4x^4+O(x^5),\qquad a_4\ne0,\qquad I_p'=p,\quad I_p(0)=0,
\]
and put
\[
F_p=\frac{y^3}{3}+y^2x^2+\frac{yxyx}{2}+yx^4+
\frac{x^6}{6}-I_p(x).
\]
\begin{theorem}\label{nrf:thm:finite}
For every cyclic formal series $h\in O_6$, the algebra
\[
A=\nrJac(F_p+h)
\]
is finite-dimensional over $\C$.

\end{theorem}

If the completed Jacobi algebra were infinite-dimensional, its associated graded algebra would be infinite-dimensional as well. In that case, the associated graded algebra would appear as a quotient of an auxiliary graded algebra having a central simple generic fiber of dimension \(25\) over the parameter field, and this quotient would necessarily inherit the same generic rank \(25\). The fourth-order relation would then force \(4\mid25\). To establish this divisibility, we let a lift of the graded central element act on the left. Right multiplication by elements of the algebra is linear with respect to this left parameter action, even when the lifted element itself is not central.

\subsection{The auxiliary graded algebra}
Assign weights \(1,2\) to \(x,y\), respectively, and define
\[
Q=y^2+yx^2+x^2y+xyx+x^4,\qquad S=yxy+y^2x+yx^3,
\qquad B=\C\langle x,y\rangle/(x^4,Q).
\]
The free polynomial identities
\[
[S,y]=y[x,Q],\qquad [S,x]=[y,Q]+x[x,Q]
\]
demonstrate that \(S\) commutes with both generators modulo \(Q=0\). Consequently, \(S\) is central in the quotient algebra \(B\) and has weighted degree five.

\begin{lemma}\label{nrf:lem:aux}
The algebra $B$ is finite over $\C[S]$, and
\[
\dim_\C B/(S)=40,\qquad B/(S-1)\cong M_5(\C).
\]
Consequently $B\otimes_{\C[S]}\C(S)$ is central simple of dimension $25$.
\end{lemma}
\begin{proof}
We compute a noncommutative Gr\"obner basis for the defining relations of \(B\), ordering words first by weighted degree and then lexicographically. The script \texttt{nr\_B5\_generic\_certificate.py} performs this completion using exact rational arithmetic, resolving every overlap and inclusion. Setting \(S=0\) yields \(40\) irreducible words, while setting \(S=1\) yields \(25\) irreducible words. Both calculations take place in the uncompleted free associative algebra over \(\Q\); extending scalars to \(\C\) gives the stated complex algebras.

For \(B/(S-1)\), the regular trace pairing is \((a,b)\mapsto\operatorname{Tr}(L_{ab})\), where \(L_{ab}\) denotes left multiplication by \(ab\) on this \(25\)-dimensional algebra. Its determinant is
\[
298023223876953125=5^{25},
\]
and the center of the specialized algebra has dimension one. In characteristic zero, nondegeneracy of the trace pairing implies that the algebra is semisimple. Because the center is one-dimensional, there is a unique simple factor, which is isomorphic to the matrix algebra \(M_5(\C)\). The file \texttt{nr\_B5\_generic\_certificate.json} lists the bases, irreducible words and exact outputs.

Lifting the homogeneous basis words of \(B/(S)\) to \(B\), induction on weighted degree shows that these lifts generate \(B\) as a module over the polynomial ring \(\C[S]\). To determine the generic fiber, adjoin an element \(z\) satisfying \(z^5=S\) to obtain the field extension \(\C(S) \subset \C(z)\). Rescaling the generators by \(x=zX\) and \(y=z^2Y\) identifies \(B \otimes_{\C[S]} \C(z)\) with the scalar extension of \(B/(S-1)\) to the field \(\C(z)\). Central simplicity and dimension \(25\) descend along this field extension.
\end{proof}

\subsection{Filtered modules with a noncentral parameter}
Let $A$ be a separated complete algebra with a decreasing positive integer filtration $F^nA$, $F^0A=A$, finite-dimensional associated graded pieces, and $F^iA F^jA\subseteq F^{i+j}A$. Thus \(\bigcap_nF^nA=0\) and \(A\cong\varprojlim A/F^nA\). For a module over \(D=\C[[t]]\) or \(\C[\sigma]\), rank means its dimension after tensoring with the fraction field of that ring.

\begin{lemma}\label{nrf:lem:rank}
Suppose $\operatorname{gr} A$ is finite over the polynomial ring $\C[\sigma]$, where
$\sigma$ is central and homogeneous of degree $e>0$. Choose a lift
$s\in F^eA$ of $\sigma$. Regard $A$ as a left $D=\C[[t]]$-module by
$t\cdot a=sa$. Then $A$ is a finite $D$-module and
\[
\operatorname{rank}_D A
=\operatorname{rank}_{\C[\sigma]}\operatorname{gr} A.
\]
Its $D$-torsion submodule is a right ideal and is finite-dimensional over
$\C$.
\end{lemma}
\begin{proof}
Lift a finite homogeneous generating family of the graded module \(\operatorname{gr}A\) over \(\C[\sigma]\) to elements of \(A\). Express any element of \(A\) by subtracting left powers of \(s\) times these lifts to eliminate the leading graded term, and repeat this division on the remainder. Because the filtration degrees strictly increase at each step, the coefficients of the lifts converge in the complete ring \(D\). This establishes that \(A\) is finitely generated as a left \(D\)-module. If \(c\) denotes the maximum filtration degree among the chosen generators, this division process yields the containment bounds
\[
 F^{en+c}A\subseteq s^nA\subseteq F^{en}A .
\]
The division process expresses an element of \(F^{en+c}A\) using left powers of \(s\) with exponent at least \(n\). A finitely generated graded \(\C[\sigma]\)-module of rank \(r\) has cumulative vector space dimension \(rN/e+O(1)\) in degrees strictly below \(N\), as follows from choosing a homogeneous basis modulo torsion. The two filtration inclusions therefore imply \(\dim_\C A/s^nA=rn+O(1)\). By the structure theorem for finitely generated modules over the discrete valuation ring \(D\), the leading coefficient of \(n\) in the dimension growth coincides with \(\operatorname{rank}_D A\).

By associativity of the algebra multiplication, right multiplication by any element of \(A\) commutes with the left action of \(s\). Consequently, the \(D\)-torsion submodule of \(A\), consisting of elements killed by some power of \(t\), is preserved by all right multiplications, making it a right ideal. Being a finitely generated torsion module over the discrete valuation ring \(D\), this submodule has finite length, and is therefore finite-dimensional over \(\C\).
\end{proof}

\begin{remark}
Left multiplication by $s$ supplies the $D$-module structure, while
right multiplication supplies the $D$-linear determinant operators.
\end{remark}

\Needspace{9\baselineskip}
\begin{lemma}\label{nrf:lem:det}
In the situation of Lemma~\ref{nrf:lem:rank}, assume the common rank is $r>0$.
Suppose $z\in F^1A$, $m\ge1$, $c\in\C^\times$, and a one-variable
formal series $v$ satisfy
\[
q(z)=cz^m(1+zv(z))=s+\epsilon,\qquad \epsilon\in F^{e+1}A.
\]
Then $m$ divides $r$.
\end{lemma}
\begin{proof}
The relation \(q(z)=s+\epsilon\) gives determinant valuation \(r\) for right multiplication by \(q(z)\), because the initial form of \(s\) is the central parameter of degree \(e\). The factorization \(q(z)=cz^m(1+zv(z))\) gives a valuation divisible by \(m\). We prove these two evaluations on the same free \(D\)-module.

\noindent\emph{(i) Right multiplication on a free module.}
Let \(T \subset A\) denote the finite-dimensional \(D\)-torsion submodule. The quotient module \(M=A/T\) is then a free module of rank \(r\) over the discrete valuation ring \(D\). We equip \(M\) with the quotient filtration \(F^n M = (F^n A + T)/T\). Because the torsion submodule \(T\) is a right ideal, right multiplication \(R_a: m \mapsto ma\) by any element \(a\in A\) descends to a well-defined endomorphism on \(M\), which is \(D\)-linear with respect to the left \(D\)-action. Centrality of the initial form of \(s\) in \(\operatorname{gr}A\) means that \([s,F^nA]\subseteq F^{n+e+1}A\) for every \(n\). This yields the inclusion
\[
(R_s-t\,\mathrm{id})(F^nM)\subseteq F^{n+e+1}M.
\]
Consequently, right multiplication by \(q(z)\) takes the form \(R_{q(z)}=t\,\mathrm{id}+E\), where the error operator \(E\) shifts the filtration by at least \(e+1\).

\par\smallskip\noindent\emph{(ii) A lattice for the filtration.}
Adjoin an element \(u\) satisfying \(t=u^e\) and work over the Laurent field \(\C((u))\). To translate the filtration shifts into valuations, consider the lattice
\[
\mathcal L=\sum_{n\ge0}u^{-n}
\bigl(F^nM\otimes_D\C[[u]]\bigr)
\subset M\otimes_D\C((u)).
\]
Set \(M_u=M\otimes_D\C[[u]]\), and let \(c_0\) be the largest weight of the homogeneous generators used in the proof of Lemma~\ref{nrf:lem:rank}. The filtration comparison in that proof gives \(F^nM\subset t^{\lfloor(n-c_0)/e\rfloor}M\) for \(n\ge c_0\). Every summand in the definition of \(\mathcal L\) is therefore contained in \(u^{-c_0-e}M_u\), while the summand at \(n=0\) is \(M_u\). The submodule \(\mathcal L\) is therefore a full \(\C[[u]]\)-lattice: it is finite free and spans \(M\otimes_D\C((u))\) over \(\C((u))\).

\par\smallskip\noindent\emph{(iii) The two determinant valuations.}
On this lattice, any operator that increases filtration degree by \(j\) maps \(\mathcal L\) into \(u^j\mathcal L\). In particular, the error operator satisfies \(E\mathcal L\subset u^{e+1}\mathcal L\), while \(t\) acts as scalar multiplication by \(u^e\). Choose a basis of \(\mathcal L\), and let \(v_t,v_u\) denote the valuations normalized by \(v_t(t)=v_u(u)=1\). For some matrix \(C\) over \(\C[[u]]\), we obtain
\[
R_{q(z)}=u^e(I+uC)\quad\text{on }\mathcal L,
\qquad v_t(\det R_{q(z)})=r.
\]
Indeed, the determinant of \(R_{q(z)}\) has \(u\)-valuation \(er\), and the valuations satisfy \(v_u=e\,v_t\) on \(\C((t))\). Because \(z \in F^1 A\), the right multiplication operator \(R_z\) maps \(\mathcal L\) into \(u\mathcal L\). The operator \(R_{1+zv(z)} = \mathrm{id} + R_z v(R_z)\) therefore reduces to the identity modulo \(u\), which gives it unit determinant in \(\C[[u]]^\times\). Taking determinants in the factorization \(q(z) = c z^m (1+zv(z))\) gives
\[
r=m\,v_t(\det R_z).
\]
Because the determinant \(\det R_z\) is a nonzero element of \(\C((t))\), its valuation \(v_t(\det R_z)\) is an integer. Thus \(m\) divides \(r\).
\end{proof}

\subsection{The rank contradiction}
\begin{proof}[Proof of Theorem~\ref{nrf:thm:finite}]
Decompose the formal tail as \(h=g(x)+h_{\mathrm{mix}}\), where \(g\in x^6\C[[x]]\) and every word in \(h_{\mathrm{mix}}\) contains at least one occurrence of \(y\). Any mixed word of ordinary word length at least six has \((1,2)\)-weighted degree at least seven. We absorb the pure part of the tail into \(\widetilde p=p-g'\), whose leading term remains \(a_4x^4\). Direct cyclic differentiation of the potential \(F_p\) yields
\[
\partial_yF_p=Q,\qquad \partial_xF_p=S+xQ-p(x).
\]
Setting the cyclic derivatives of \(F_p + h\) to zero in the completed Jacobi algebra \(A\), the two defining relations become
\[
Q\in F^5A,\qquad
\widetilde p(x)=S+
\bigl(\partial_x h_{\mathrm{mix}}-x\partial_yh_{\mathrm{mix}}\bigr),
\]
where the parenthesized discrepancy has weighted degree at least six. The leading terms of these two relations with respect to the \((1,2)\)-grading are \(x^4\) and \(Q\). There is therefore a canonical surjection of graded algebras
\[
B\twoheadrightarrow\operatorname{gr} A.
\]
By Lemma~\ref{nrf:lem:aux}, the associated graded algebra \(\operatorname{gr} A\) is a finitely generated module over the image of the central parameter ring \(\C[S]\).

If \(A\) were infinite-dimensional over \(\C\), its associated graded algebra \(\operatorname{gr} A\) would be infinite-dimensional as well. As a finitely generated module over \(\C[S]\), it would then necessarily have a nonzero generic fiber, because a finitely generated torsion module over a one-variable polynomial ring is finite-dimensional. This generic fiber is a nonzero quotient of the central simple algebra \(B \otimes_{\C[S]} \C(S)\) of dimension \(25\) from Lemma~\ref{nrf:lem:aux}. Because a central simple algebra admits no proper nonzero two-sided ideals, this quotient map is an isomorphism, which forces the generic rank of \(\operatorname{gr} A\) over \(\C[S]\) to be exactly \(25\).

We now apply Lemma~\ref{nrf:lem:rank} to \(s=S(x,y)\in F^5A\), whose initial form is the image \(\sigma\) of the graded central element \(S\). Lemma~\ref{nrf:lem:det} then applies with parameters \(e=5\), \(z=x\), \(q=\widetilde p\), and leading exponent \(m=4\). This forces \(4\) to divide \(25\), which is impossible. We conclude that \(A\) is finite-dimensional over \(\C\).
\end{proof}

\subsection{Polynomial truncation and length-five realization}
\label{subsec:nrf-realisation}

\begin{proposition}\label{nrf:prop:geom}
Every potential in Theorem~\ref{nrf:thm:finite} has a geometric Jacobi
algebra, realized by
a smooth length-five $cE_8$ flop. Equivalently, the alternating-quartic
stratum with nonzero pure fifth coefficient is realized, including all
formal tails.
\end{proposition}

Theorem~\ref{nrf:thm:finite} establishes that the Jacobson radical is nilpotent. The following lemma translates this nilpotence exponent into an effective polynomial truncation bound.
\begin{lemma}\label{nrf:lem:truncate}
If $\nrJac(f)$ has radical nilpotence exponent $L$, so that $J^L=0$, then replacing $f$ by
its ordinary truncation through any degree $N\ge L+1$ leaves its closed
Jacobi ideal unchanged.
\end{lemma}
\begin{proof}
Let \(I\) and \(I_N\) denote the closed two-sided ideals generated by the cyclic derivatives of \(f\) and its degree-\(N\) polynomial truncation \(f_{\le N}\), respectively. Their generators agree modulo the maximal ideal power \(\mathfrak m^N\), and by assumption \(\mathfrak m^L\subseteq I\). Because \(N \ge L+1\), we have the containment
\[
\mathfrak m^L\subseteq I_N+\mathfrak m^{L+1}.
\]
Multiplying this containment successively by powers of \(\mathfrak m\) and substituting back recursively yields \(\mathfrak m^L\subseteq I_N+\mathfrak m^{L+k}\) for arbitrarily large integers \(k\). Because the ideal \(I_N\) is closed in the \(\mathfrak m\)-adic topology, taking the intersection over all \(k\) yields \(\mathfrak m^L\subseteq I_N\). The differences between the cyclic derivatives of \(f\) and those of \(f_{\le N}\) belong to \(\mathfrak m^N \subseteq \mathfrak m^L \subseteq I \cap I_N\), which proves the equality of derivative ideals \(I=I_N\).
\end{proof}

\begin{proof}[Proof of Proposition~\ref{nrf:prop:geom}]
Applying the coordinate change \(X=y+x^2\) and \(Y=x\) transforms the potential into a series with five-jet
\[
\frac{X^3}{3}+\frac{(XY)^2}{2}-XY^4-\frac{a_4}{5}Y^5.
\]
At the level of five-jets, we perform a shear \(Y\mapsto Y+dX\). We restore the standard quartic form using a quadratic correction in \(X\) and clear newly introduced cyclic \(XX\) terms. This transformation alters the coefficient \(B\) of the necklace \(XY^4\) to \(B+5Cd\), where \(C\ne0\) denotes the coefficient of the pure fifth power \(Y^5\). We choose \(d\) to cancel this term completely. A subsequent quadratic correction in \(Y\) removes the remaining \(XYXY^2\) term. Rescaling the generators and multiplying the potential by a nonzero scalar yields \(X^3+(XY)^2+Y^5+O_6\), where all higher-order terms remain in \(O_6\).

By Theorem~\ref{nrf:thm:finite}, the Jacobson radical of the Jacobi algebra has a finite nilpotence exponent \(L\). We truncate the potential through a degree \(N\) satisfying \(N\ge\max(5,L+1)\). By Lemma~\ref{nrf:lem:truncate}, this polynomial truncation preserves the closed Jacobi ideal, while the degree bound preserves the normalized five-jet.

We then apply Wemyss's gluing theorem \cite[Sections~1.5--1.6, Eq.~(1.3) and Theorem~1.7]{WemyssTaming}. The fifth family listed in Theorem~1.6 of that paper is \(X^3+(XY)^2+Y^5+O_6\). The realization statement following Eq.~(1.3) identifies the Jacobi algebra of a listed polynomial potential with the noncommutative deformation algebra of the associated two-chart curve neighborhood. Theorem~1.7 of \cite{WemyssTaming} proves that contracting the central zero section produces a projective crepant contraction from a regular threefold source, forming a smooth length-five flop. Localizing and completing at the closed image of the curve and pulling back the regular source produces the required complete local contraction of type \(cE_8\).
\end{proof}

\begin{remark}
For an explicit instance of the gluing construction, consider the two-series subfamily \(G=F_p+g(x)+\ell(x)y\) with \(g\in x^6\C[[x]]\) and \(\ell\in x^5\C[[x]]\), and set \(\widetilde p=p-g'\). Take coordinates \((a,u,v)\) and \((b,U,V)\) on the two charts, where \(a,b\) are the usual affine coordinates on \(\mathbb P^1\). The gluing construction gives the following transition functions. With formal series coefficients, these coordinate transitions are interpreted on the completed charts by
\[
b=a^{-1},\quad U=a^{-1}u,\quad
V=a^3v+a^{-1}u^2+u^3+au^4+a^2(u^5-\widetilde p(u))+a\ell(u).
\]
This coordinate transition has Jacobian determinant \(-1\). Its zero section is a curve \(\mathbb P^1\) with normal bundle \(\mathcal O(1)\oplus\mathcal O(-3)\). By \cite[Corollary 6.4]{BrownWemyssDerived}, the deformation algebra of this curve neighborhood is isomorphic to \(\nrJac(G)\). Truncating the series yields polynomial gluing data for this subfamily.
\end{remark}

\begin{remark}
The subsequent classification of Brown and Wemyss provides an alternative deduction of the geometric realization step. In cyclic necklace notation, let \(\mathcal N_{2,2}=X^2Y^2+(XY)^2/2\). Eliminating the bad quartic term from \(X^3+\mathcal N_{2,2}+Y^5\) leaves a nonzero alternating quartic coefficient and a nonzero pure fifth coefficient. The coordinate transformations in Corollary~\ref{cor:nr-five-jet} put this five-jet and the five-jet used in Proposition~\ref{nrf:prop:geom} into the same normal form. Reversing the changes for \(X^3+\mathcal N_{2,2}+Y^5\) gives the coordinate identification established in \cite[Lemma~8.1]{BrownWemyssClass}. Theorem~8.3 of that paper carries every formal tail of this five-jet to an isomorphic Jacobi algebra belonging to the fifth polynomial family in Corollary~8.14. Finally, Theorem~9.2 of \cite{BrownWemyssClass}, implication \((2)\Rightarrow(3)\), realizes that family geometrically.
\end{remark}

\section{Growth exclusions}
\label{sec:gr-growth-exclusions}
We set $A_f=\nrJac(f)$ and write $J=\operatorname{rad} A_f$ for the Jacobson radical of $A_f$. Any continuous coordinate change with continuous inverse preserves the radical filtration, which ensures that the quotient dimensions $\dim_\C A_f/J^N$ are invariant under such automorphisms.

\begin{theorem}[Zero-quartic exclusion]\label{nrg:thm:zero}
For every cyclic formal potential $f=x^3+O_5$ and every $N\ge1$,
\[
 \dim_\C A_f/J^N\ \ge\
 Z(N):=\sum_{k=0}^{\lfloor(N-1)/3\rfloor}\binom{N-3k+1}{2}
 =\frac{N^3}{18}+O(N^2).
\]
In particular $\nrJdim A_f\ge3$.
\end{theorem}

\begin{theorem}[Both degenerate sextic axes]\label{nrg:thm:boundary}
Let
\[
 f=x^3+y^4+b(xy^2)^2+c(xy)^3+O_7,\qquad bc=0,
\]
with an arbitrary complete cyclic tail. Then $\nrJdim A_f\ge3$.
More precisely, after the normalization with weights \(\operatorname{wt}(x,y)=(4,3)\) in Lemma~\ref{nrg:lem:normalise}, the following bounds hold for every positive integer $N$, and hence
also hold for the original algebra:
\begin{align*}
 b=0:\quad &\dim_\C A_f/J^{3N}\ge 3\binom{N+2}{3},\\
 c=0:\quad &\dim_\C A_f/J^{2N}\ge
 2\sum_{k=0}^{\lfloor(N-1)/2\rfloor}\binom{N-2k+1}{2}.
\end{align*}
The intersection $b=c=0$ satisfies both bounds.
\end{theorem}

The proof of Theorem~\ref{nrg:thm:zero} begins with an associated graded calculation on a Zariski open set of fifth jets, followed by specialization at each finite cutoff. For the degenerate sextic axes in Theorem~\ref{nrg:thm:boundary}, splitting a pure power produces a complete quiver whose Jacobi algebra yields lower bounds for the original radical quotients.

\subsection{The associated graded algebra on an open set of fifth jets}
\label{sec:gr-pbw}

The next proposition computes the algebra associated to the radical filtration, \(\operatorname{gr}_J A=\bigoplus_{n\ge0}J^n/J^{n+1}\), for a normalized fifth jet and every cyclic tail of order at least six. For a graded algebra \(B=\bigoplus_{n\ge0}B_n\), we use the Hilbert series \(H_B(t)=\sum_{n\ge0}(\dim B_n)t^n\). Subsection~\ref{sec:gr-zero-growth} identifies the open set represented by this jet and extends the dimension bound by specialization.

\begin{proposition}\label{nrg:prop:genericzero}
For every cyclic $h\in O_6$, let
\[
 f=\frac{x^3}{3}+xyxy^2+\frac{y^5}{5}+h.
\]
Then
\[
 \operatorname{gr}_J\nrJac(f)\cong
 \frac{\C\langle x,y\rangle}
 {(x^2,\ y^4+xy^2x+yxyx+xyxy)}.
\]
In particular, $\operatorname{gr}_J\nrJac(f)$ has Hilbert series $1/((1-t)^2(1-t^3))$.
\end{proposition}

To prove Proposition~\ref{nrg:prop:genericzero}, we determine the homogeneous relations and their syzygies over \(\Q\), and then lift these syzygies after extending scalars to \(\C\).
Let $R=\Q\langle x,y\rangle$ carry the standard grading by word length. We define the homogeneous elements
\[
 U=yxy^2+y^2xy,\qquad
 R_0=x^2,\qquad
 R_1=y^4+xy^2x+yxyx+xyxy.
\]
Let $M=RpR\oplus RqR$ be the free $R$-bimodule on formal symbols $p,q$ of weights $2,4$. The $R$-bimodule homomorphism $d_0:M\to R$ substitutes $R_0$ for $p$ and $R_1$ for $q$. An element of $\ker d_0$ is a syzygy: a linear combination of two-sided relation multiples that evaluates to zero.

\begin{lemma}\label{nrg:lem:certificate}
The kernel of $d_0$ is generated as an $R$-bimodule by
\begin{align*}
 s_1&=xp-px,\\
 \sigma&=x(yq-qy)+(yq-qy)x-Up+pU,
\end{align*}
and the following Peiffer identities, which compare the two orders of applying a pair of relations
\[
 R_i a\theta_j-\theta_i aR_j,
 \qquad i,j\in\{0,1\},\quad a\in R,
 \quad\theta_0=p,\ \theta_1=q.
\]
Moreover
\[
 H_{R/(R_0,R_1)}(t)=\frac1{(1-t)^2(1-t^3)}.
\]
\end{lemma}
\begin{proof}
There are two assertions to check. The critical identities generate the relation kernel, and the irreducible words determine the Hilbert series. We use the degree-lexicographic order with $x<y$. A homogeneous Gr\"obner basis for the two-sided ideal is
\begin{align*}
 &x^2,\\
 &y^4+yxyx+xy^2x+xyxy,\\
 &y^3xy^2x+y^3xyxy-yxyxy^3-xy^2xy^3,\\
 &y^3xyx+y^2xy^2x+y^2xyxy-yxyxy^2-xy^2xy^2-xyxy^3,\\
 &y^2xyx+yxy^2x-xy^2xy-xyxy^2.
\end{align*}
The fourth basis element is redundant. All five appear in the certificate. The Python verifier checks that all five elements lie in the ideal of the two original relations and resolves all seventeen overlaps and inclusions. Substituting the expressions in the certificate back into the free relation bimodule gives eleven nonzero critical syzygies. Their weights and multiplicities are
\[
 \begin{array}{c|rrrrrr}
 \text{weight}&3&6&7&8&9&10\\\hline
 \text{number}&1&1&1&4&3&1.
 \end{array}
\]
The JSON certificate expresses these eleven critical syzygies using thirty-one explicit rational multiples of $s_1,\sigma$, and Peiffer identities. The verifier checks these identities in the free relation bimodule over the polynomial free algebra, using exact rational coefficients.

The reduction from the full kernel to these finite identities proceeds by induction on the maximal leading word. In any polynomial syzygy among the relations, we compare two occurrences of relations matching at the maximal leading word. When the occurrences overlap or include one another, their difference gives a critical syzygy corresponding to an overlap or inclusion ambiguity. When the occurrences are disjoint, their difference forms a Peiffer identity. Subtracting the verified resolution of each ambiguity strictly lowers the maximal leading word, ensuring that the reduction terminates after finitely many steps. Replacing each auxiliary Gr\"obner generator with its polynomial expression in $R_0, R_1$ from the certificate makes its defining relation vanish identically, while expanding its Peiffer identities produces sums of Peiffer identities for the original generators $R_0, R_1$. The induction and the verified identities show that $s_1, \sigma$, and the Peiffer identities generate $\ker d_0$ as an $R$-bimodule.

\par\smallskip\noindent\emph{The Hilbert series.}
The leading monomials of the Gröbner basis yield the minimal forbidden words
\[
 xx,\qquad yyyy,\qquad yyxyx,\qquad yyyxyyx.
\]
To count words avoiding these four factors, a prefix automaton keeps track of the longest suffix of the current word that is a proper prefix of a forbidden word. There are ten such states, including the empty word. A transition \(i\to j\) appends one letter without creating a forbidden factor. The certificate provides exact polynomial numerators $n_i(t)$ with common denominator $D(t)=(1+t)(1-t)^2(1-t^3)$ satisfying the transition equations
\[
 n_i(t)=D(t)+t\sum_{i\to j}n_j(t).
\]
Here \(n_i(t)/D(t)\) counts the allowed continuations from state \(i\), including the empty continuation. At the state corresponding to the empty word, the numerator is $1+t$. Because the linear system defined by these automaton equations admits a unique formal power series solution, these verified polynomial identities prove the asserted rational Hilbert series in every degree.
\end{proof}

Tensor each homogeneous component of \(M\xrightarrow{d_0}R\) with \(\C\) over \(\Q\). In each degree this is a map between finite-dimensional rational vector spaces. Flatness of scalar extension preserves its kernel, so the displayed syzygies and their bimodule multiples generate the relation kernel over \(\C\) as well. From this point on, set \(R=\C\langle x,y\rangle\) and \(\widehat R=\prod_{n\ge0}R_n\), its completion by word length. The initial ideal of a closed ideal is generated by the lowest-degree nonzero homogeneous parts of all its elements. The following lemma proves that no additional initial relations arise when the syzygies of the given initial relations lift exactly. Its proof works modulo one finite power of the augmentation ideal at a time.

\begin{lemma}\label{nrg:lem:lift}
Let $r_i$ be filtered formal relations with homogeneous initial relations
$R_i$. If every homogeneous syzygy of the $R_i$ is a sum of homogeneous
syzygies admitting exact syzygy lifts for the $r_i$ with the same initial
terms, then
\[
 \operatorname{gr}\bigl(\widehat R/\overline{(r_i)}\bigr)\cong R/(R_i).
\]
\end{lemma}
\begin{proof}
Let an element of the closed relation ideal have order $n$, with nonzero initial form in degree $n$. Modulo order $n+1$, we choose a representation of this element as a sum of relation multiples. If the lowest-weight terms of this representation cancel, that cancelling component forms a homogeneous syzygy of the initial relations $R_i$. By hypothesis, every such homogeneous syzygy is a sum of homogeneous syzygies that admit exact lifts to syzygies of the filtered relations $r_i$ with the same initial forms. Subtracting these exact lifts leaves the represented element unchanged while strictly increasing the minimal weight appearing in its representation. The least weight strictly increases and cannot pass $n$, so finitely many subtractions express the weight-$n$ initial form as a combination of the $R_i$. The initial ideal of $\overline{(r_i)}$ is therefore contained in $(R_i)$. Conversely, the initial ideal contains each $R_i$ as the initial form of $r_i$. This argument applies to the topologically closed ideal because determining the initial form of an element of order $n$ requires only a representation valid through degree $n$.
\end{proof}

\begin{proof}[Proof of Proposition~\ref{nrg:prop:genericzero}]
By Lemma~\ref{nrg:lem:lift}, it suffices to lift the generators of the homogeneous relation kernel supplied by Lemma~\ref{nrg:lem:certificate}. The two cyclic derivatives of $f$ are
\[
 r_x=x^2+U+\partial_xh,\qquad
 r_y=R_1+\partial_yh.
\]
Assign weights $2,4$ to the formal relation symbols $p,q$, and let $d$ evaluate these symbols at $r_x,r_y$, respectively. The identity \([x,\partial_xf]+[y,\partial_yf]=0\) gives the exact syzygy
\[
 S_1=[x,p]+[y,q],\qquad dS_1=[x,r_x]+[y,r_y]=0.
\]
The initial term is $s_1$ of weight three. Using \(S_1\), define a second element of the completed relation bimodule
\[
 S_2=xS_1+S_1x-r_xp+pr_x.
\]
Indeed, direct evaluation shows $dS_2=0$. Cancelling its weight-four terms gives
\[
 S_2=\sigma-(\partial_xh)p+p(\partial_xh).
\]
Since $\partial_xh$ has order at least five, the last two terms have weight at least seven. Thus $S_2$ lifts the weight-six syzygy $\sigma$. Replacing $R_i$ by $r_i$ lifts every Peiffer identity. The scalar-extension argument above Lemma~\ref{nrg:lem:lift} shows that these lifts cover every homogeneous syzygy over \(\C\). The lifting lemma therefore identifies the initial ideal with $(R_0, R_1)$ and proves the asserted graded isomorphism.

\end{proof}

\subsection{Specialization to all fifth jets}
\label{sec:gr-zero-growth}
\begin{proof}[Proof of Theorem~\ref{nrg:thm:zero}]
We first identify the Zariski open set of fifth jets to which Proposition~\ref{nrg:prop:genericzero} applies. The relation matrices at a finite cutoff will then extend the resulting Hilbert-function bound to every fifth jet.
Scaling the potential allows us to assume that its cubic part is \(x^3/3\). In degree five, an order-three coordinate correction of $x$ removes every necklace containing $xx$. The remaining fifth jet has the form
\[
 \frac{x^3}{3}+Axyxy^2+Bxy^4+Cy^5.
\]
When $C\ne0$, we apply the linear shear $y\mapsto y+tx$ and remove the cyclic $xx$ terms generated by this substitution. The surviving coefficients become
\[
 A'=A+2Bt+5Ct^2,\qquad B'=B+5Ct,\qquad C'=C.
\]
The shear leaves the discriminant $B^2-5AC$ invariant. On the open set defined by
\[
 C\ne0,\qquad B^2-5AC\ne0,
\]
we take $t=-B/(5C)$. This choice kills the $xy^4$ term, since $B'=0$, while leaving $A'=-(B^2-5AC)/(5C)\ne0$. Rescaling $x, y$ and the overall potential recovers the exact fifth jet of Proposition~\ref{nrg:prop:genericzero}. No quartic terms are introduced, and the remaining $O_6$ tail is arbitrary. The coordinate changes induce continuous isomorphisms of the completed Jacobi algebras. Multiplying the potential by a nonzero scalar leaves its derivative ideal unchanged. Both operations preserve the radical quotient dimensions, so Proposition~\ref{nrg:prop:genericzero} gives the same truncated dimensions throughout this open set.

\par\smallskip\noindent\emph{Specialization at a finite cutoff.}
Fix \(N\ge1\). We allow the coefficients to vary through degree \(\max(N,5)\), so that the conditions on the fifth jet remain defined even for small \(N\). Modulo \(\mathfrak m^N\), the ambient algebra is finite-dimensional, and the two-sided multiples of the cyclic derivatives define a relation matrix whose entries are polynomials in the potential coefficients. On the dense open set of fifth jets, Proposition~\ref{nrg:prop:genericzero} proves that the quotient dimension $\dim_\C A_f/J^N$ is constant, which means that the relation matrix has constant rank \(r\). Every minor of size \(r+1\) vanishes on this open set and hence vanishes identically throughout the parameter space. At any specialized potential, the relation matrix can therefore only have rank at most \(r\). Because the quotient dimension is the ambient dimension minus the rank of the relation matrix, every specialization has quotient dimension at least the generic value. The cyclic derivatives of potential terms beyond degree \(N\) lie in \(\mathfrak m^N\), so this establishes the lower bound for the complete formal potential.

Finally, writing \([t^n]\) for extraction of the coefficient of \(t^n\), the Hilbert series gives
\[
 [t^n]\frac1{(1-t)^2(1-t^3)}
 =\#\{(a,b,c)\in\mathbb Z_{\ge0}^3:a+b+3c=n\}.
\]
Summing over $n<N$ gives $Z(N)$. The leading quadratic terms of the binomial coefficients sum to $N^3/18+O(N^2)$, which completes the proof of Theorem~\ref{nrg:thm:zero}.
\end{proof}

\subsection{Two complete-quiver growth bounds}
\label{sec:gr-quiver-growth}
Let $Q_q$ denote the quiver with $q$ vertices and one arrow for each ordered pair of vertices, including a loop at each vertex. Over a field \(k\), we work in the path algebra of \(Q_q\) completed with respect to the ideal \(\mathfrak a\) generated by its arrows. A quiver potential is a possibly infinite linear combination of cyclic paths, modulo cyclic rotation; its Jacobi algebra is the completed path algebra modulo its cyclic arrow derivatives. We also write \(\mathfrak a\) for the image of the arrow ideal in this quotient. Dimensions below are taken over \(k\). The two degenerate sextic axes in Theorem~\ref{nrg:thm:boundary} correspond to the pairs $(q,m)=(3,3)$ and $(2,4)$, where the potential has order at least $m$ in the arrow ideal.

\begin{proposition}\label{nrg:prop:quiver}
Over any characteristic-zero extension of $\C$, the following hold.
For every potential $W\in\mathfrak a^3$ on $Q_3$,
\[
 \dim \nrJac(W)/\mathfrak a^N\ge3\binom{N+2}{3}.
\]
For every potential $W\in\mathfrak a^4$ on $Q_2$,
\[
 \dim \nrJac(W)/\mathfrak a^N\ge
 2\sum_{k=0}^{\lfloor(N-1)/2\rfloor}\binom{N-2k+1}{2}.
\]
The potentials may be formal and need not be homogeneous.
\end{proposition}

The estimates follow by varying the potential coefficients about two explicit homogeneous potentials. At each finite cutoff, vanishing in cohomological degree minus one on the special fiber gives flatness over the parameter ring.

We use the Ginzburg differential graded (dg) algebra of the quiver potential. For each arrow \(a:i\to j\), adjoin an opposite arrow \(a^*:j\to i\) of cohomological degree \(-1\), and at each vertex adjoin a loop \(t_i\) of degree \(-2\). The original arrows have degree zero. With \(e_i\) the vertex idempotents, the differential is determined by
\[
 d(a)=0,\qquad d(a^*)=\partial_aW,\qquad
 d(t_i)=e_i\Bigl(\sum_a[a,a^*]\Bigr)e_i,
\]
and the graded Leibniz rule. Its zeroth cohomology is the completed Jacobi algebra. In addition, positive internal weights yield finite truncations. At the two homogeneous potentials constructed below, the differential preserves internal weight exactly.

\begin{lemma}\label{nrg:lem:flat}
Let $(T,\mathfrak n,k)$ be a Noetherian local ring and let $P$ be a
nonpositive complex of flat $T$-modules. Suppose $H^0(P)$ is finite and
$H^{-1}(P\otimes_Tk)=0$. Then $H^0(P)$ is finite free over $T$.
\end{lemma}
\begin{proof}
It suffices to show that \(\operatorname{Tor}_1^T(H^0(P),k)=0\). Put $K=\operatorname{im}(P^{-1}\to P^0)$ and $B=H^0(P) = P^0/K$. Given an element of $K\otimes_T k$ that vanishes in $P^0\otimes_T k$, lift it to $P^{-1}\otimes_T k$. The image of this lift in $P^0\otimes_T k$ is zero, so the lift is a cycle. By the hypothesis that $H^{-1}(P\otimes_Tk)=0$, this cycle is a boundary coming from $P^{-2}\otimes_T k$, which implies that its image in $K\otimes_T k$ vanishes. Hence the canonical map $K\otimes_T k\to P^0\otimes_T k$ is injective. From the short exact sequence $0\to K\to P^0\to B\to 0$, the flatness of $P^0$ now gives $\operatorname{Tor}_1^T(B,k)=0$.

Take a minimal surjection from a finite free $T$-module onto $B$. Its kernel is finitely generated by Noetherianity and has zero reduction modulo $\mathfrak n$ by the Tor vanishing and minimality. Nakayama's lemma then implies that the kernel is zero, which proves the claim.
\end{proof}

\begin{lemma}\label{nrg:lem:seeds}
There are homogeneous potentials $W_3$ of degree three on $Q_3$ and $W_2$
of degree four on $Q_2$ whose Ginzburg differential graded algebras have
zero negative cohomology. Their Jacobi algebras have total Hilbert series
\[
 H_3(t)=\frac{3}{(1-t)^3},\qquad
 H_2(t)=\frac{2}{(1-t)^2(1-t^2)}.
\]
For either homogeneous potential, of degree $m$, assign internal weight one to every arrow;
opposite generators have weight
$m-1$ and the degree-minus-two vertex generators have weight $m$.
\end{lemma}
\begin{proof}
We perform the initial calculations within the graded path algebras, after which completion takes the direct product of the homogeneous components. Because the differential for each homogeneous potential preserves the internal weight grading, these calculations simultaneously determine the completed algebras in each internal weight and at every finite internal-weight cutoff.

\par\smallskip\noindent\emph{(i) The cubic potential on \(Q_3\).}
For $q=3$, let the cyclic group of order three act on the polynomial ring $\C[X,Y,Z]$ with characters $0,1,2$. The McKay quiver has the group characters as vertices; a generator of character \(j\) gives an arrow from character \(i\) to character \(i+j\). For the characters \(0,1,2\) of \(X,Y,Z\), this is precisely \(Q_3\). Viewing \(XYZ-XZY\) as a cyclic potential in the free noncommuting generators, its McKay lift is a quiver potential whose uncompleted Jacobi algebra is \(\C[X,Y,Z]\rtimes C_3\). Here \(C_3\) is the cyclic group of order three, and the skew group product is \((ag)(bh)=a\,g(b)\,gh\). Explicitly, taking vertex indices modulo three, let $a_i:i\to i$ denote the loops, $b_i:i\to i+1$ the forward arrows, and $c_i:i\to i-1$ the backward arrows. The cubic potential is
\[
 W_3=\sum_i(a_ib_ic_{i+1}-a_ic_ib_{i-1}).
\]
The cyclic derivatives give the three commutation relations with their vertex labels. The graded Ginzburg differential graded algebra is the cobar algebra of the exterior Koszul coalgebra on three generators, crossed with $C_3$. Here the cobar algebra is the tensor dg algebra on the shifted reduced coalgebra; the exterior degree-\(j\) generators have cohomological degree \(1-j\), and the differential is induced by the coproduct. The Koszul resolution identifies its augmentation to the skew polynomial algebra as a quasi-isomorphism. Because the volume element is invariant, the degree-minus-two generators are vertex loops. The ordered monomials \(X^iY^jZ^k\), together with the three group elements, give the first Hilbert series.

\par\smallskip\noindent\emph{(ii) The quartic potential on \(Q_2\).}
For $q=2$, consider the order-two group $C_2$ acting by the involution $X\mapsto X$ and $Y\mapsto-Y$. The McKay lift of the invariant potential $\frac12[X,Y]^2$ is a quartic potential on the two-vertex quiver $Q_2$. Before completion, the resulting Jacobi algebra is
\[
 U(\mathfrak h_3)\rtimes C_2,\qquad
 [X,Y]=Z,\quad [X,Z]=[Y,Z]=0,
\]
where \(\mathfrak h_3\) is the three-dimensional Heisenberg Lie algebra with the displayed brackets, \(U(\mathfrak h_3)\) is its universal enveloping algebra, and $\deg X=\deg Y=1$, $\deg Z=2$. The involution sends $(X,Y,Z)$ to $(X,-Y,-Z)$. The Poincar\'e--Birkhoff--Witt (PBW) theorem gives the basis \(X^iY^jZ^k\). With the two group elements and the weights \(1,1,2\), this yields the second Hilbert series. To prove the required vanishing of negative cohomology, we realize this algebra through a three-generator potential and then eliminate a contractible pair. Regard \(Z\) as a third free generator and use the stabilized potential
\[
 \widetilde W=Z[X,Y]-\tfrac12Z^2.
\]
The associated Ginzburg differential graded algebra is the cobar construction of the Chevalley--Eilenberg coalgebra of $\mathfrak h_3$, crossed with $C_2$. This coalgebra is the exterior coalgebra on the Lie algebra, with differential induced by its Lie bracket. The cobar augmentation sends the degree-zero generators to the corresponding Lie algebra elements in \(U(\mathfrak h_3)\) and the negative-degree generators to zero. This map is an equivariant quasi-isomorphism. Indeed the PBW filtration reduces it to the exterior Koszul calculation for the abelian three-dimensional Lie algebra. Since the differential preserves the positive internal weights, the conclusion holds separately in each weight.

The remaining step is to return to the quartic Ginzburg algebra. Denote the opposite generators of the stabilized algebra by $p,q,r$, so that
\[
 dp=[Y,Z],\quad dq=[Z,X],\quad dr=[X,Y]-Z,
 \quad d\zeta=[X,p]+[Y,q]+[Z,r].
\]
Set
\[
 Z'=Z-[X,Y],\quad p'=p+[Y,r],\quad q'=q+[r,X],\quad \zeta'=\zeta+r^2.
\]
Then $dr=-Z'$, $dZ'=0$, and
\[
 dp'=[Y,[X,Y]],\quad dq'=[[X,Y],X],\quad
 d\zeta'=[X,p']+[Y,q'].
\]
Thus the stabilized differential graded algebra is the quartic Ginzburg algebra with the pair \((r,Z')\) freely adjoined. The equations \(dr=-Z'\) and \(dZ'=0\) make this pair a two-term contractible complex. Projection onto the quartic differential graded algebra is a quasi-isomorphism. In each internal weight, the words involving the pair form tensor complexes with a contractible factor. The coordinate change respects both the weights and the involution, proving the required vanishing.
\end{proof}

\begin{proof}[Proof of Proposition~\ref{nrg:prop:quiver}]
For a fixed arrow cutoff we first obtain a free family near the chosen homogeneous potential. The polynomial relation matrix then gives the dimension bound everywhere. Fix \(N\) and \((q,m)\in\{(3,3), (2,4)\}\). Only potential terms in degrees \(m,\ldots,N\) affect \(\nrJac(W)/\mathfrak a^N\). Vary their cyclic coefficients independently about \(W_q\) over \(T = \C[[e_1,\ldots,e_M]]\). Give the arrows, opposite arrows and vertex loops of the relative Ginzburg differential graded algebra over \(T\) internal weights \(1, m-1, m\), respectively. The differential never lowers internal weight. We can therefore discard weights at least \(N\) in every cohomological degree, leaving a complex of finite free \(T\)-modules.

At the homogeneous potential \(W_q\), the differential preserves internal weight exactly. Its truncated complex is the direct sum of the full weight complexes below \(N\), so Lemma~\ref{nrg:lem:seeds} gives vanishing in cohomological degree minus one. In cohomological degree zero, the cohomology of the relative complex is the truncated Jacobi algebra \(\nrJac(W_T)/\mathfrak a^N\). The complex has finite free terms and its zeroth cohomology is finite. Lemma~\ref{nrg:lem:flat} therefore makes this cohomology free over \(T\), with rank given by its special fiber.

The differential into cohomological degree zero is a matrix of polynomials in the potential coefficients. If its generic rank is \(r\), every \((r+1)\times(r+1)\) minor vanishes as a formal series at the coefficient point of \(W_q\) and hence vanishes identically. At every coefficient point the rank is at most \(r\), so the quotient dimension is at least the dimension at \(W_q\). Translating the coefficient coordinates by those of \(W_q\) covers the full coefficient space. Summing the coefficients of the two Hilbert series in Lemma~\ref{nrg:lem:seeds} gives the asserted bounds. The minor identities survive scalar extension, and terms of degree greater than \(N\) do not affect this quotient.
\end{proof}

\subsection{Growth transfer by partial morsification}
\label{sec:gr-morsification}
We deform \(v^{q+1}/(q+1)\) to \(v^{q+1}/(q+1)-s^qv\), whose derivative \(v^q-s^q\) has \(q\) distinct roots after \(s\) is inverted. This is the partial morsification used here: only the critical points in the \(v\)-variable are separated. Each root gives a vertex, and the blocks of the remaining generator \(u\) give the arrows between vertices. We compare the original algebra with the generic algebra at all \(q\) critical points.

\Needspace{19\baselineskip}
\begin{lemma}\label{nrg:lem:transfer}
Let
\[
 f=P(u)+\frac{v^{q+1}}{q+1}+g(u,v),
\]
where $P\in u^m\C[[u]]$ and every word of $g$ contains at least $m$
occurrences of $u$, with $m\ge3$. Suppose that, for each fixed number of
$u$'s, only finitely many words of $g$ occur; the same conclusion holds if
the substitutions below converge in the corresponding adic topology.
Put $A=\nrJac(f)$, let $I$ be its closed two-sided ideal generated by $u$, and put $J=\operatorname{rad} A$.
For each $N$ there is a finite $\C[[s]]$-module $M_N$ whose special fiber is
$A/I^N$ and whose generic fiber is
\[
 \nrJac(W)/\mathfrak a^N
\]
for a potential $W\in\mathfrak a^m$ on $Q_q$ over $K=\C((s))$.
In particular,
\[
 \dim_\C A/J^{qN}\ \ge\ \dim_\C A/I^N
 \ \ge\ \dim_K\nrJac(W)/\mathfrak a^N.
\]

\end{lemma}
\begin{proof}
We construct a finite module over \(D=\C[[s]]\). Its special fiber is \(M_N/sM_N\), and its generic fiber is \(M_N\otimes_D\C((s))\). We identify the latter at all \(q\) roots of \(v^q-s^q\), so that none of the critical points is lost by completion at a single root.

\noindent\emph{(i) The finite module over \(D\).}
Take \(F_s=f-s^qv\) over \(D=\C[[s]]\) and let \(M_N\) be its relative Jacobi quotient modulo \(I^N\), completed in the total \((s,u,v)\)-adic topology. Its special fiber is \(A/I^N\). Since \(A/I=\C[[v]]/(v^q)\), we have \(J^{qN}\subset I^N\). Put \(K=\C((s))\) and \(\lambda_i=s\zeta_q^i\), where \(\zeta_q\) is a primitive \(q\)-th root of unity.

To keep the full generic support, we first work over the parameter ring $D$. The relative relation \(v^q-s^q=-\partial_vg\in I^m\) gives
\[
                     (v^q-s^q)^{\lceil N/m\rceil}=0
                         \qquad\pmod{I^N}.
\]
Under the finite-word hypothesis the truncated derivative relations are polynomials, since only terms with at most \(N\) occurrences of \(u\) contribute. We form the polynomial quotient modulo \(I^N\) over \(\C[[s]]\). Its reduction modulo \(I\) is the finite free algebra \(\C[[s]][v]/(v^q-s^q)\). Each subquotient \(I^j/I^{j+1}\) is spanned by products of \(j\) copies of \(u\) separated by representatives of a basis of this algebra. The polynomial quotient is therefore finite over the discrete valuation ring \(\C[[s]]\). Modulo \(s\) its augmentation ideal is nilpotent, which implies that the total-adic and \(s\)-adic topologies are cofinal. The quotient is already complete and equals \(M_N\).

For convergent substitutions we take the remainder of each string of \(v\)'s between successive \(u\)'s modulo the displayed monic polynomial in \(v\). The lower coefficients lie in \((s)\), so the remainders converge \(s\)-adically and establish the same finite-module property.

\par\smallskip\noindent\emph{(ii) The completed algebra at all \(q\) roots.}
Finiteness of $M_N$ over the parameter discrete valuation ring $D$ bounds the dimension of the generic fiber by that of the special fiber. To identify the generic fiber, we invert \(s\). Modulo \(I\) we obtain \(K^q\). The powers of the factors \(v-\lambda_i\) are pairwise coprime for the $q$ distinct roots, so the Chinese remainder theorem gives idempotents \(e_1,\ldots,e_q\) that separate the \(q\) roots. They are polynomials in \(v\), satisfy \(e_ie_j=\delta_{ij}e_i\), and sum to one. These idempotents commute with \(v\), although they need not commute with \(u\).

We use these spectral idempotents before imposing the Jacobi relations. We complete the free algebra over \(K\) at the kernel of the map \(K\langle u,v\rangle\to K^q\) sending \(u\) to zero and \(v\) to \((\lambda_1,\ldots,\lambda_q)\). Successive lifting in the \((v^q-s^q)\)-adic topology gives the orthogonal idempotents inside the closed commutative subalgebra generated by \(v\). Relative to this copy of \(K[v]/(v^q-s^q) \cong \prod_{i=1}^q K\), the completed free algebra has generators
\[
 u_{ij}=e_iue_j\quad(1\le i,j\le q),\qquad
 t_i=e_i(v-\lambda_i)e_i\quad(1\le i\le q).
\]
These generate the completed path algebra of the quiver obtained from \(Q_q\) by adjoining one additional loop \(t_i\) at each vertex. The \(q^2\) arrows \(u_{ij}\) belong to \(Q_q\) itself. The inverse substitution sends \(u\) to \(\sum_{i,j} u_{ij}\) and \(v\) to \(\sum_i(\lambda_i e_i+t_i)\). The spectral projectors of this last element recover the vertex idempotents.

\par\smallskip\noindent\emph{(iii) Elimination of the quadratic loops.}
The expansion of the pure \(v\)-potential \(\frac{v^{q+1}}{q+1} - s^q v\) at each root \(\lambda_i\) has nondegenerate quadratic part in the loop \(t_i\), because
\[
 (v^q-s^q)'|_{v=\lambda_i}=q\lambda_i^{q-1}\ne0.
\]
Every other term in \(f\) contains at least \(m\) copies of \(u\). We complete the squares successively in the loops \(t_i\) while fixing the arrows \(u_{ij}\), following the noncommutative splitting lemma. The critical loop values satisfy \(t_i\in(u_{ij})^m\). Substituting them therefore produces a quiver potential \(W\) of arrow order at least \(m\), that is, \(W \in \mathfrak a^m\).

Before eliminating the quadratic loops, differentiation with respect to the arrows and loops gives the matrix blocks of \(\partial_u F_s\) and the diagonal blocks of \(\partial_v F_s\). The splitting imposes the loop equations and leaves the derivatives of \(W\) in the surviving arrows. After the arrow derivative relations impose \(\partial_uF_s=0\), the identity \([u,\partial_uF_s]+[v,\partial_vF_s]=0\) gives \([v,\partial_vF_s]=0\). On any off-diagonal block with \(i\ne j\), the commutator operator
\[
 z\longmapsto(\lambda_i+t_i)z-z(\lambda_j+t_j)
\]
is invertible. Its constant term is the scalar \(\lambda_i-\lambda_j = s(\zeta_q^i - \zeta_q^j) \ne 0\), and the remaining operator is topologically nilpotent. Thus the off-diagonal blocks of the derivative vanish in the quotient too. The resulting isomorphism preserves the ideal generated by the arrows \(u_{ij}\).

\par\smallskip\noindent\emph{(iv) The dimension comparison.}
Quotienting by the \(N\)-th power of the arrow ideal yields a finite-dimensional \(K\)-algebra. The substitutions identify the generic fiber of \(M_N\) with the arrow-ideal quotient \(\nrJac(W)/\mathfrak a^N\). Indeed the ideal generated by \(u\) and \(v^q-s^q\) is nilpotent modulo \(I\), and \(I^N=0\), so completion at the full semisimple quotient preserves this finite fiber. The finiteness of $M_N$ over the discrete valuation ring $D$ establishes the dimension inequality $\dim_\C A/I^N \ge \dim_K \nrJac(W)/\mathfrak a^N$. Combining this inequality with the inclusion \(J^{qN}\subset I^N\) completes the proof of both bounds.
\end{proof}

\begin{remark}
The finite-word hypothesis holds for the normalized potentials below. In a cyclic word avoiding $xx$ and $yyy$, the number of $y$'s is at most twice the number of $x$'s, and the number of $x$'s is at most the number of $y$'s. A bound on either letter count therefore bounds the word length.
\end{remark}

\subsection{Application to the sextic axes}
\label{sec:gr-sextic}
\begin{lemma}\label{nrg:lem:normalise}
Give $x,y$ weights $4,3$. Every $x^3+y^4+O_5$ is right-equivalent
to a potential whose higher necklaces avoid cyclic $xx$
and cyclic $yyy$. Below weight $24$ its only possible higher necklaces are
\[
 a\,xyxy^2\quad(\operatorname{wt}17),\qquad
 b(xy^2)^2\quad(\operatorname{wt}20),\qquad
 c(xy)^3\quad(\operatorname{wt}21).
\]
For $x^3+y^4+b(xy^2)^2+c(xy)^3+O_7$, the normalization leaves $b,c$
unchanged and creates no weight-$17$ term.
\end{lemma}
\begin{proof}
Lemma~\ref{lem:nr-y4-reduction} gives the convergent weighted normalization. For the displayed sextic jet there are no earlier bad terms to alter its coefficients. The terms of ordinary order $O_7$ in weights $21,22,23$ are all bad necklaces. Removing them changes no good necklace of those weights, so the normalized remainder still has ordinary order at least seven.
\end{proof}

\begin{proof}[Proof of Theorem~\ref{nrg:thm:boundary}]
Lemma~\ref{nrg:lem:normalise} preserves the two sextic coefficients and removes the bad necklaces. For these normalized tails, a bound on either letter count bounds the total word length. This verifies the finite-word hypothesis of Lemma~\ref{nrg:lem:transfer}. We choose \(u=x\) on the axis \(b=0\) and \(u=y\) on the axis \(c=0\); the required letter count is different in the two cases.

First suppose $b=0$. Every mixed word in the normalized tail contains at least three occurrences of $x$. Indeed, with only two $x$'s a good necklace has at most four $y$'s, hence length at most six. The only such sextic is $(xy^2)^2$, absent because $b=0$, and the weight-seventeen necklace $xyxy^2$ is absent as well. After rescaling the coordinates we have
\[
 f=\frac{x^3}{3}+\frac{y^4}{4}+g,\qquad
 \text{where every word of }g\text{ contains at least three }x\text{'s}.
\]
We apply Lemma~\ref{nrg:lem:transfer} with $u=x$, $v=y$, and parameters $(q,m)=(3,3)$. The uniform quiver bound of Proposition~\ref{nrg:prop:quiver} yields
\[
 \dim_\C A/J^{3N}\ge3\binom{N+2}{3}.
\]

If $c=0$, every mixed term of the normalized tail contains at least four occurrences of $y$. A good necklace with at most three $y$'s has at most three $x$'s and thus length at most six. The only possible sextic is $(xy)^3$, absent because $c=0$. We now take $u=y$ and split the pure power in $v=x$. Applying Lemma~\ref{nrg:lem:transfer} with parameters $(q,m)=(2,4)$ together with Proposition~\ref{nrg:prop:quiver} yields
\[
 \dim_\C A/J^{2N}\ge
 2\sum_{k=0}^{\lfloor(N-1)/2\rfloor}\binom{N-2k+1}{2}.
\]
Both bounds grow cubically in $N$, so $\nrJdim A_f \ge 3$ on either axis. At their intersection $b=c=0$ both bounds apply.
\end{proof}

\section{Proof of the realization theorem}
\label{sec:gr-assembly}

\begin{proof}[Proof of Theorem~\ref{thm:nrf-numerical-realisation}]
By Proposition~\ref{prop:nr-full-reduction} it suffices to realize the potentials \(x^3+O_4\) satisfying the numerical growth bound. We apply the quartic classification of Proposition~\ref{prop:nr-quartic-orbits}.

If \(\Delta\ne0\), Theorem~\ref{thm:nr-e6-complete} constructs a length-three arc for each potential in this stratum and separates the small contractions from the divisorial contractions.

For an alternating quartic, Theorem~\ref{thm:nr-weighted-growth} demonstrates that any potential with $\nrJdim A_f \le 1$ must have a nonzero pure fifth coefficient. We normalize the five-jet using Corollary~\ref{cor:nr-five-jet}. Reversing the coordinate comparison in the proof of Proposition~\ref{nrf:prop:geom}, which uses \(X=y+x^2\), \(Y=x\), puts this normalized potential into the form \(F_p+h\) with \(p(u)=u^4\) and \(h\in O_6\). The coordinate map and its inverse preserve ordinary order. Theorem~\ref{nrf:thm:finite} proves finite-dimensionality, and Proposition~\ref{nrf:prop:geom} realizes the algebra geometrically as a length-five flop.

For the pure \(y^4\) quartic orbit, Lemma~\ref{nrg:lem:normalise} isolates the surviving higher coefficients in weights \(17,20,21\). If the weight-\(17\) coefficient is nonzero, Theorem~\ref{thm:nr-length-four} realizes its Jacobi algebra by a length-four contraction. If the weight-\(17\) coefficient vanishes, Theorem~\ref{nrg:thm:boundary} excludes both degenerate sextic axes $b=0$ and $c=0$, leaving exactly the open pair of nonzero coefficients realized in length six by Theorem~\ref{thm:nr-length-six}. The remaining quartic orbit consisting of potentials with zero quartic jet is excluded by the cubic growth bound of Theorem~\ref{nrg:thm:zero}.

\par\smallskip\noindent\textbf{Generic surface types.}
Varying the hyperplane section in the fixed threefold deforms its distinguished partial-resolution surface while preserving the scheme-theoretic fiber. The generic Du Val diagram is a Dynkin subgraph of the distinguished diagram with the same multiplicity at the marked vertex. In the length-three \(E_6\) family no proper Dynkin subgraph supports multiplicity three. Thus the generic singularity type is \(cE_6\) for both the small contractions and the divisorial contractions. For the small primitive contractions of lengths four, five, and six, the invariant-theoretic criteria of Katz and Morrison \cite[Section~8]{KatzMorrison} establish the generic singularity types as \(cE_7,cE_8,cE_8\), respectively.

Under $\nrJdim\le1$ the surviving Type E branches are therefore
\[
\begin{array}{lll}
\toprule
\text{normalized branch}&\text{construction}&\text{generic surface type}\\
\midrule
\Delta\ne0&\text{marked length }3&cE_6\\
(xy)^2\text{ quartic, nonzero pure fifth}&\text{marked length }5&cE_8\\
y^4\text{ quartic, }\alpha_{17}\ne0&\text{marked length }4&cE_7\\
y^4\text{ quartic, }\alpha_{17}=0,\ \alpha_{20}\alpha_{21}\ne0
  &\text{marked length }6&cE_8\\
\bottomrule
\end{array}
\]
The growth bounds exclude every other branch.

In lengths three, four, and six, Lemma~\ref{lem:nr-tilting-base-change} identifies the completed Jacobi algebra with the basic contraction algebra of the geometric model. For length five, Lemma~\ref{nrf:lem:truncate} identifies the closed Jacobi ideal with that of the polynomial truncation used in Proposition~\ref{nrf:prop:geom}. All continuous coordinate changes and scalar rescalings preserve the completed Jacobi algebra up to isomorphism.

Together with the Type A and Type D cases established in Proposition~\ref{prop:nr-full-reduction}, these constructions satisfy Definition~\ref{def:nr-geometric} for every algebra in the theorem. The finite-dimensional algebras are realized by flopping contractions and the growth-one algebras by divisor-to-curve contractions.
\end{proof}

\clearpage

\section{Computational certificates}
\label{sec:gr-certificates}
The programs and exact certificates are supplied as arXiv ancillary files in the folder \nolinkurl{geometric_realization_computational_code}. The separate computational ZIP archive expands to the same folder. In the complete arXiv source archive, this folder is under \texttt{anc/}. From the directory containing the folder, run
\begin{verbatim}
cd geometric_realization_computational_code
python3 tools/numerical_realisation_certificates/verify_cyclic_reductions.py
python3 tools/numerical_realisation_certificates/verify_zero_quartic.py
python3 tools/numerical_realisation_certificates/nr_B5_generic_certificate.py \
    --output output/nr_B5_rerun.json
\end{verbatim}
The programs use exact arithmetic and require Python 3.10 or later. The first two use the standard library. Before running the third, install the versions of SymPy and mpmath listed in \texttt{requirements-certificates.txt} at the top of the extracted code folder. The output directory is created automatically, and the generated JSON can be compared with the supplied length-five certificate.

For Lemma~\ref{lem:nr-e6-tangent}, \nolinkurl{verify_cyclic_reductions.py} reduces each of the ten overlap ambiguities in both reduction orders. For Theorem~\ref{thm:nr-length-six} it forms the ordinary and cyclic relation matrices of \(B_6\), computes their exact ranks over the rationals, and checks the seven proposed cocenter generators.

The script \nolinkurl{verify_zero_quartic.py} reads the adjacent \texttt{zero\_quartic\_certificate.json}. It checks the five original-ideal identities, the seventeen overlap and inclusion ambiguities, and the thirty-one summands representing the eleven critical syzygies in Lemma~\ref{nrg:lem:certificate}. It also verifies the polynomial identities giving the rational Hilbert series of the prefix automaton. Lemma~\ref{nrg:lem:lift} passes from the homogeneous relations to formal relations.

The script \nolinkurl{nr_B5_generic_certificate.py} recomputes the quotients \(B/(S)\) and \(B/(S-1)\) in Lemma~\ref{nrf:lem:aux}. The supplied \nolinkurl{nr_B5_generic_certificate.json} contains the monomial basis words, the vector space dimensions \(40\) and \(25\), and the regular trace determinant and center dimension at \(S=1\). These verified values supply the finite inputs to the filtered argument of Theorem~\ref{nrf:thm:finite}.

\section*{Acknowledgements}
\addcontentsline{toc}{section}{Acknowledgements}
I made extensive use of large language models, including ChatGPT, and the mathematical exploration agents Danus and Rethlas throughout this work. These tools assisted with mathematical exploration and exposition, suggested proof strategies and intermediate arguments, and helped check and reorganize calculations. I formulated the research questions, developed the strategies, and selected and revised the arguments presented here. I have made every effort to verify these arguments against the cited literature and accompanying exact certificates, and I take full responsibility for the manuscript and for making it public.

Neither the existence nor the priority date of this manuscript should discourage a wholly independent human proof of these results. If a subsequent work obtains them through an argument developed wholly independently by humans, I claim no special priority or credit for the proof ideas beyond standard bibliographic acknowledgement of prior or concurrent results.

If later work makes substantive use of the results, constructions, or proof strategies in this manuscript, ordinary scholarly citation is appropriate regardless of the tools used. In particular, I would appreciate explicit attribution when closely related arguments appear in AI-assisted work, since the origin of ideas generated through related language models can otherwise be difficult to trace.

\end{document}